\documentclass[12pt]{article}
\usepackage{maple2e}
\usepackage{a4}
\usepackage{amssymb}
\newcommand{\N} {{\mathbb N}}
\newcommand{\Z} {{\mathbb Z}}
\newcommand{\Q} {{\mathbb Q}}
\newcommand{\hypergeom}[5]{\mbox{$
_#1 F_#2\left.
\!\!
\left(
\!\!\!\!
\begin{array}{c}
\multicolumn{1}{c}{\begin{array}{c}
#3
\end{array}}\\[1mm]
\multicolumn{1}{c}{\begin{array}{c}
#4
            \end{array}}\end{array}
\!\!\!\!
\right| \displaystyle{#5}\right)
$}
}
\newcommand{\qfac}[2]{\left(#1;q\right)_{#2}}
\newcommand{\qphihyp}[5]{{}_{#1}\phi_{#2}\left.\left(\begin{array}{c}
        #3 \\ #4 \end{array}\right|q,#5\right)}
\newcommand{\qbinomial}[2]{\left[\begin{array}{c}\!\!#1\!\!\\
\!\!#2\!\!\end{array}\right]_{q}}

\title{Software for the Algorithmic Work with Orthogonal Polynomials and 
Special Functions}
\author{Wolfram Koepf \\
Hochschule f\"{u}r Technik, Wirtschaft und Kultur Leipzig, Germany\\[3mm]
Plenary Talk at the IWOP 98, Madrid}
\date{June 29--30, 1998}

\begin{document}
\maketitle

\section{Basics}
Modern computer algebra systems use heu\-ristics and algorithms for
the fast computation with mathematical formulas.

General purpose computer algebra systems like
{\sl Axiom} \cite{Axiom}, {\sl Derive} \cite{Derive}, {\sl Macsyma}
\cite{Macsyma}, {\sl Maple} \cite{Maple}, {\sl Mathematica} \cite{Mathematica}
or {\sl Reduce} \cite{Reduce} are for example great with integrals.
Even a small system like Derive computes {\sl all} explicitly given
integrals of Bronshtein and Semedyayev's integral table \cite{BS}.
But how do such computations work?

To begin with I would like to 
give examples of important mathematical concepts and methods that are
available in general purpose computer algebra systems.
For demonstration purposes I use the Maple V.5 system.

\subsection{Linear Algebra}

One of the main topics of any computer algebra system is linear algebra.
Linear algebra algorithms are used throughout Mathematics; we will
see examples in connection with orthogonal polynomials and special functions
later. 

With Maple, we can compute the solution of a linear system of equations:

\begin{maplegroup}
\begin{mapleinput}
\mapleinline{active}{1d}{solve(\{x+2*a*y+3*z=4,5*x+6*y+7*z=8,9*x+10*y+11*z=12\},}{%
}
\mapleinline{active}{1d}{\hspace{0mm}\{x,y,z\});}{%
}
\end{mapleinput}

\mapleresult
\begin{maplelatex}
\[
\{z={\displaystyle \frac {3}{2}} , \,x={\displaystyle \frac {-1}{
2}} , \,y=0\}
\]
\end{maplelatex}

\end{maplegroup}
\noindent
even if parameters are involved. For this purpose, Maple uses a Gauss
type algorithm.

Note, that the above system is linear only if $a$ is considered constant.
If we consider $a$ as a variable, then a nonlinear system has to be solved:

\begin{maplegroup}
\begin{mapleinput}
\mapleinline{active}{1d}{solve(\{x+2*a*y+3*z=4,5*x+6*y+7*z=8,9*x+10*y+11*z=12\},}{%
}
\mapleinline{active}{1d}{\hspace{0mm}\{a,x,y,z\});}{%
}
\end{mapleinput}

\mapleresult
\begin{maplelatex}
\begin{eqnarray*}
\lefteqn{\{z={\displaystyle \frac {3}{2}} , \,x={\displaystyle 
\frac {-1}{2}} , \,y=0, \,a=a\}, \{} \\
 & & z= - {\displaystyle \frac {1}{2}} \,y + {\displaystyle 
\frac {3}{2}} , \,a=1, \,x= - {\displaystyle \frac {1}{2}} \,y - 
{\displaystyle \frac {1}{2}} , \,y=y  \}
\;.
\end{eqnarray*}
\end{maplelatex}

\end{maplegroup}
\noindent
In a forthcoming section,
we give more details on nonlinear systems of equations.

Maple has a large linear algebra library:

\begin{maplegroup}
\begin{mapleinput}
\mapleinline{active}{1d}{with(linalg);}{%
}
\end{mapleinput}

\mapleresult
\begin{maplettyout}
Warning, new definition for norm
Warning, new definition for trace
\end{maplettyout}

\begin{maplelatex}
\begin{eqnarray*}
\lefteqn{[\mathit{BlockDiagonal}, \,\mathit{GramSchmidt}, \,
\mathit{JordanBlock}, \,\mathit{LUdecomp}, \,\mathit{QRdecomp}, }
 \\
 & & \mathit{Wronskian}, \,\mathit{addcol}, \,\mathit{addrow}, \,
\mathit{adj}, \,\mathit{adjoint}, \,\mathit{angle}, \,\mathit{
augment}, \,\mathit{backsub},  \\
 & & \mathit{band}, \,\mathit{basis}, \,\mathit{bezout}, \,
\mathit{blockmatrix}, \,\mathit{charmat}, \,\mathit{charpoly}, \,
\mathit{cholesky}, \,\mathit{col},  \\
 & & \mathit{coldim}, \,\mathit{colspace}, \,\mathit{colspan}, \,
\mathit{companion}, \,\mathit{concat}, \,\mathit{cond}, \,
\mathit{copyinto},  \\
 & & \mathit{crossprod}, \,\mathit{curl}, \,\mathit{definite}, \,
\mathit{delcols}, \,\mathit{delrows}, \,\mathit{det}, \,\mathit{
diag}, \,\mathit{diverge},  \\
 & & \mathit{dotprod}, \,\mathit{eigenvals}, \,\mathit{
eigenvalues}, \,\mathit{eigenvectors}, \,\mathit{eigenvects},  \\
 & & \mathit{entermatrix}, \,\mathit{equal}, \,\mathit{
exponential}, \,\mathit{extend}, \,\mathit{ffgausselim}, \,
\mathit{fibonacci},  \\
 & & \mathit{forwardsub}, \,\mathit{frobenius}, \,\mathit{
gausselim}, \,\mathit{gaussjord}, \,\mathit{geneqns}, \,\mathit{
genmatrix},  \\
 & & \mathit{grad}, \,\mathit{hadamard}, \,\mathit{hermite}, \,
\mathit{hessian}, \,\mathit{hilbert}, \,\mathit{htranspose}, \,
\mathit{ihermite},  \\
 & & \mathit{indexfunc}, \,\mathit{innerprod}, \,\mathit{intbasis
}, \,\mathit{inverse}, \,\mathit{ismith}, \,\mathit{issimilar}, 
\,\mathit{iszero},  \\
 & & \mathit{jacobian}, \,\mathit{jordan}, \,\mathit{kernel}, \,
\mathit{laplacian}, \,\mathit{leastsqrs}, \,\mathit{linsolve}, \,
\mathit{matadd},  \\
 & & \mathit{matrix}, \,\mathit{minor}, \,\mathit{minpoly}, \,
\mathit{mulcol}, \,\mathit{mulrow}, \,\mathit{multiply}, \,
\mathit{norm}, \,\mathit{normalize},  \\
 & & \mathit{nullspace}, \,\mathit{orthog}, \,\mathit{permanent}
, \,\mathit{pivot}, \,\mathit{potential}, \,\mathit{randmatrix}, 
 \\
 & & \mathit{randvector}, \,\mathit{rank}, \,\mathit{ratform}, \,
\mathit{row}, \,\mathit{rowdim}, \,\mathit{rowspace}, \,\mathit{
rowspan}, \,\mathit{rref},  \\
 & & \mathit{scalarmul}, \,\mathit{singularvals}, \,\mathit{smith
}, \,\mathit{stackmatrix}, \,\mathit{submatrix}, \,\mathit{
subvector},  \\
 & & \mathit{sumbasis}, \,\mathit{swapcol}, \,\mathit{swaprow}, 
\,\mathit{sylvester}, \,\mathit{toeplitz}, \,\mathit{trace}, \,
\mathit{transpose},  \\
 & & \mathit{vandermonde}, \,\mathit{vecpotent}, \,\mathit{
vectdim}, \,\mathit{vector}, \,\mathit{wronskian}]
\end{eqnarray*}
\end{maplelatex}

\end{maplegroup}
\noindent
You can see which procedures are available now. As an example, we
compute the determinant of the matrix
\[
\left(
\begin{array}{ccc}
1&2a&3\\
4&5&6\\
7&8&9
\end{array}
\right)
\]
by

\begin{maplegroup}
\begin{mapleinput}
\mapleinline{active}{1d}{det([[1,2*a,3],[5,6,7],[9,10,11]]);}{%
}
\end{mapleinput}

\mapleresult
\begin{maplelatex}
\[
 - 16 + 16\,a
\]
\end{maplelatex}

\end{maplegroup}
\noindent
and the eigenvalues and eigenvectors for $a=1$:

\begin{maplegroup}
\begin{mapleinput}
\mapleinline{active}{1d}{eigenvalues([[1,2,3],[5,6,7],[9,10,11]]);}{%
}
\end{mapleinput}

\mapleresult
\begin{maplelatex}
\[
0, \,9 + \sqrt{105}, \,9 - \sqrt{105}
\]
\end{maplelatex}

\end{maplegroup}
\begin{maplegroup}
\begin{mapleinput}
\mapleinline{active}{1d}{eigenvectors([[1,2,3],[5,6,7],[9,10,11]]);}{%
}
\end{mapleinput}

\mapleresult
\begin{maplelatex}
\begin{eqnarray*}
\lefteqn{[0, \,1, \,\{[1, \,-2, \,1]\}], [9 + \sqrt{105}, \,1, }
 \\
 & & \{ \left[  \! {\displaystyle \frac {11}{2}}  - 
{\displaystyle \frac {1}{2}} \,\sqrt{105}, \,1, \, - 
{\displaystyle \frac {7}{2}}  + {\displaystyle \frac {1}{2}} \,
\sqrt{105} \!  \right] \}], [ \\
 & & 9 - \sqrt{105}, \,1,  \\
 & & \{ \left[  \! {\displaystyle \frac {11}{2}}  + 
{\displaystyle \frac {1}{2}} \,\sqrt{105}, \,1, \, - 
{\displaystyle \frac {7}{2}}  - {\displaystyle \frac {1}{2}} \,
\sqrt{105} \!  \right] \}]
\;.
\end{eqnarray*}

\end{maplelatex}

\end{maplegroup}
\noindent
Later we will show how important an efficient implementation of
linear algebra can be.

\subsection{Polynomial Artithmetic}

A second major topic of computer algebra is polynomial arithmetic. 

$P$ defines a polynomial

\begin{maplegroup}
\begin{mapleinput}
\mapleinline{active}{1d}{P:=(1-x)*sum(x^k,k=0..9);}{%
}
\end{mapleinput}

\mapleresult
\begin{maplelatex}
\[
P := (1 - x)\,(1 + x + x^{2} + x^{3} + x^{4} + x^{5} + x^{6} + x
^{7} + x^{8} + x^{9})
\]
\end{maplelatex}

\end{maplegroup}
\noindent
which is expanded by

\begin{maplegroup}
\begin{mapleinput}
\mapleinline{active}{1d}{expand(P);}{%
}
\end{mapleinput}

\mapleresult
\begin{maplelatex}
\[
1 - x^{10}
\;.
\]
\end{maplelatex}

\end{maplegroup}
\noindent
$Q$ is a second polynomial

\begin{maplegroup}
\begin{mapleinput}
\mapleinline{active}{1d}{Q:=1-x^5;}{%
}
\end{mapleinput}

\mapleresult
\begin{maplelatex}
\[
Q := 1 - x^{5}
\]
\end{maplelatex}

\end{maplegroup}
\noindent
and {\tt normal} cancels common factors of the ratio $P/Q$:

\begin{maplegroup}
\begin{mapleinput}
\mapleinline{active}{1d}{normal(P/Q);}{%
}
\end{mapleinput}

\mapleresult
\begin{maplelatex}
\[
x^{5} + 1
\]
\end{maplelatex}

\end{maplegroup}
\noindent
by an application of the Euclidean algorithm.

A highlight of computer algebra is {\sl rational factorization} since the
underlying algorithms (factorization modulo a prime and Hensel lifting,
or the triple $L$ algorithm, see e.g.\ 
\cite{GCL}) are not suitable for hand computations.

For example, the polynomial $P$ can be factored over $\Q$ by the command

\begin{maplegroup}
\begin{mapleinput}
\mapleinline{active}{1d}{factor(P);}{%
}
\end{mapleinput}

\mapleresult
\begin{maplelatex}
\[
 - ( - 1 + x)\,(x + 1)\,(x^{4} + x^{3} + x^{2} + x + 1)\,(x^{4}
 - x^{3} + x^{2} - x + 1)
\]
\end{maplelatex}

\end{maplegroup}
\noindent
and the following is a rational factorization of $1-x^{105}$: 

\begin{maplegroup}
\begin{mapleinput}
\mapleinline{active}{1d}{factor(1-x^105);}{%
}
\end{mapleinput}

\mapleresult
\begin{maplelatex}
\begin{eqnarray*}
\lefteqn{ - ( - 1 + x)\,(x^{6} + x^{5} + x^{4} + x^{3} + x^{2} + 
x + 1)\,(x^{4} + x^{3} + x^{2} + x + 1)(1 - x + x^{5}} \\
 & & \mbox{} - x^{6} + x^{7} - x^{8} + x^{10} - x^{11} + x^{12}
 - x^{13} + x^{14} - x^{16} + x^{17} - x^{18} + x^{19} \\
 & & \mbox{} - x^{23} + x^{24})(x^{2} + x + 1)\,(1 - x + x^{3} - 
x^{4} + x^{6} - x^{8} + x^{9} - x^{11} + x^{12}) \\
 & & (1 - x + x^{3} - x^{4} + x^{5} - x^{7} + x^{8})(1 + x + x^{2
} - x^{5} - x^{6} - 2\,x^{7} - x^{8} - x^{9} \\
 & & \mbox{} - x^{24} + x^{12} + x^{13} + x^{14} + x^{16} + x^{17
} + x^{15} + x^{48} - x^{20} - x^{22} - x^{26} - x^{28} \\
 & & \mbox{} + x^{31} + x^{32} + x^{33} + x^{34} + x^{35} + x^{36
} - x^{39} - x^{40} - 2\,x^{41} - x^{42} - x^{43} + x^{46} \\
 & & \mbox{} + x^{47})
\end{eqnarray*}
\end{maplelatex}

\end{maplegroup}
\noindent
Note that 105 is the
smallest exponent such that the rational factorization of $1-x^n$ contains
coefficients different from 0 or $\pm 1$.

Next, we define a multivariate polynomial

\begin{maplegroup}
\begin{mapleinput}
\mapleinline{active}{1d}{Product(x^(2*k-1)-y^k/k^2,k=1..7);}{%
}
\end{mapleinput}

\mapleresult
\begin{maplelatex}
\[
{\displaystyle \prod _{k=1}^{7}} \,(x^{(2\,k - 1)} - 
{\displaystyle \frac {y^{k}}{k^{2}}} )
\]
\end{maplelatex}

\end{maplegroup}
\noindent
whose expanded form is a huge expression:

\begin{maplegroup}
\begin{mapleinput}

\mapleinline{active}{1d}{term:=expand(product(x^(2*k-1)-y^k/k^2,k=1..7));
}{%
}
\end{mapleinput}

\mapleresult
\begin{maplelatex}
\begin{eqnarray*}
\lefteqn{\mathit{term} :=  - {\displaystyle \frac {1}{49}} \,x^{
36}\,y^{7} - {\displaystyle \frac {1}{25}} \,x^{40}\,y^{5} + 
{\displaystyle \frac {1}{1225}} \,x^{27}\,y^{12} + 
{\displaystyle \frac {421}{176400}} \,x^{29}\,y^{11}} \\
 & & \mbox{} - {\displaystyle \frac {1}{44100}} \,x^{16}\,y^{18}
 - {\displaystyle \frac {1}{16}} \,x^{42}\,y^{4} + 
{\displaystyle \frac {113}{28224}} \,x^{31}\,y^{10} - 
{\displaystyle \frac {1}{28224}} \,x^{18}\,y^{17} \\
 & & \mbox{} + {\displaystyle \frac {16969}{1587600}} \,x^{33}\,y
^{9} + x^{49} - {\displaystyle \frac {181}{1587600}} \,x^{20}\,y
^{16} - {\displaystyle \frac {71}{235200}} \,x^{22}\,y^{15} \\
 & & \mbox{} + {\displaystyle \frac {1}{705600}} \,x^{9}\,y^{22}
 - {\displaystyle \frac {1}{9}} \,x^{44}\,y^{3} + {\displaystyle 
\frac {5609}{176400}} \,x^{35}\,y^{8} - {\displaystyle \frac {47
}{45360}} \,x^{24}\,y^{14} \\
 & & \mbox{} + {\displaystyle \frac {1}{396900}} \,x^{11}\,y^{21}
 + {\displaystyle \frac {161}{3600}} \,x^{37}\,y^{7} - 
{\displaystyle \frac {5099}{3175200}} \,x^{26}\,y^{13} + 
{\displaystyle \frac {61}{6350400}} \,x^{13}\,y^{20} \\
 & & \mbox{} - {\displaystyle \frac {2069}{564480}} \,x^{28}\,y^{
12} + {\displaystyle \frac {1}{26880}} \,x^{15}\,y^{19} + 
{\displaystyle \frac {13}{181440}} \,x^{17}\,y^{18} - 
{\displaystyle \frac {1}{6350400}} \,x^{4}\,y^{25} \\
 & & \mbox{} - {\displaystyle \frac {1}{4}} \,x^{46}\,y^{2} + 
{\displaystyle \frac {89}{1600}} \,x^{39}\,y^{6} - 
{\displaystyle \frac {17147}{3175200}} \,x^{30}\,y^{11} + 
{\displaystyle \frac {559}{3628800}} \,x^{19}\,y^{17} \\
 & & \mbox{} - {\displaystyle \frac {1}{2822400}} \,x^{6}\,y^{24}
 + {\displaystyle \frac {13}{144}} \,x^{41}\,y^{5} - 
{\displaystyle \frac {18733}{1587600}} \,x^{32}\,y^{10} + 
{\displaystyle \frac {167}{453600}} \,x^{21}\,y^{16} \\
 & & \mbox{} - {\displaystyle \frac {13}{6350400}} \,x^{8}\,y^{23
} - {\displaystyle \frac {21}{1600}} \,x^{34}\,y^{9} + 
{\displaystyle \frac {2629}{5080320}} \,x^{23}\,y^{15} \\
 & & \mbox{} - {\displaystyle \frac {89}{25401600}} \,x^{10}\,y^{
22} + {\displaystyle \frac {1091}{1270080}} \,x^{25}\,y^{14} - 
{\displaystyle \frac {1}{90720}} \,x^{12}\,y^{21} \\
 & & \mbox{} - {\displaystyle \frac {209}{12700800}} \,x^{14}\,y
^{20} + {\displaystyle \frac {1}{25401600}} \,x\,y^{27} - y\,x^{
48} + {\displaystyle \frac {1}{9}} \,y^{4}\,x^{43} \\
 & & \mbox{} - {\displaystyle \frac {61}{3600}} \,y^{8}\,x^{36}
 + {\displaystyle \frac {401}{313600}} \,y^{13}\,x^{27} - 
{\displaystyle \frac {23}{635040}} \,y^{19}\,x^{16} + 
{\displaystyle \frac {1}{6350400}} \,y^{26}\,x^{3} \\
 & & \mbox{} + {\displaystyle \frac {1}{4}} \,y^{3}\,x^{45} - 
{\displaystyle \frac {1}{64}} \,y^{7}\,x^{38} + {\displaystyle 
\frac {181}{129600}} \,y^{12}\,x^{29} - {\displaystyle \frac {113
}{2822400}} \,y^{18}\,x^{18} \\
 & & \mbox{} + {\displaystyle \frac {1}{2822400}} \,y^{25}\,x^{5}
 - {\displaystyle \frac {1}{36}} \,y^{6}\,x^{40} + 
{\displaystyle \frac {1}{900}} \,y^{11}\,x^{31} - {\displaystyle 
\frac {421}{6350400}} \,y^{17}\,x^{20} \\
 & & \mbox{} + {\displaystyle \frac {1}{1587600}} \,y^{24}\,x^{7}
 + {\displaystyle \frac {1}{576}} \,y^{10}\,x^{33} - 
{\displaystyle \frac {1}{20736}} \,y^{16}\,x^{22} + 
{\displaystyle \frac {1}{1016064}} \,y^{23}\,x^{9} \\
 & & \mbox{} - {\displaystyle \frac {1}{36}} \,x^{38}\,y^{6} + 
{\displaystyle \frac {1}{1764}} \,x^{25}\,y^{13} - 
{\displaystyle \frac {1}{14400}} \,y^{15}\,x^{24} + 
{\displaystyle \frac {1}{705600}} \,y^{22}\,x^{11} \\
 & & \mbox{} + {\displaystyle \frac {1}{518400}} \,y^{21}\,x^{13}
 - {\displaystyle \frac {1}{25401600}} \,y^{28}
\end{eqnarray*}
\end{maplelatex}

\end{maplegroup}
\noindent
which is a polynomial of degree $49$ w.r.t.\ $x$ and of degree $28$
w.r.t.\ $y$.

It is beautiful (and will turn out to be essential in the sequel)
that computer algebra systems
have no problems to factorize such expressions over the rationals in
reasonable time:

\begin{maplegroup}
\begin{mapleinput}
\mapleinline{active}{1d}{factor(term);}{%
}
\end{mapleinput}

\mapleresult
\begin{maplelatex}
\begin{eqnarray*}
\lefteqn{{\displaystyle \frac {1}{25401600}} (x - y)\,( - y^{2}
 + 4\,x^{3})\,( - y^{3} + 9\,x^{5})\,( - y^{5} + 25\,x^{9})\,( - 
y^{6} + 36\,x^{11})} \\
 & & (49\,x^{13} - y^{7})\,( - y^{4} + 16\,x^{7})
\mbox{\hspace{171pt}}
\end{eqnarray*}
\end{maplelatex}

\end{maplegroup}

\subsection{Polynomial Systems}
\label{subsection:Polynomial Systems}

We come back to the problem of nonlinear systems of equations. 
Whereas in the linear case, Gauss elimination works, Buchberger's
algorithm is an extension to the multivariate case. It
constitutes---given a certain term order---an elimination scheme
to find a normal form for a given polynomial system, which can be used
to find the general solution of a nonlinear system.

We consider the following system of equations:

\begin{maplegroup}
\begin{mapleinput}
\mapleinline{active}{1d}{LIST:=\{9*B*A+4*d-6*c*d=0,}{%
}
\mapleinline{active}{1d}{\hspace{1.7cm}-9*a*b+9*B*A=0,}{%
}
\mapleinline{active}{1d}{\hspace{1.7cm}-18*B*A+12*d-12*c*d+4*d^2=0,}{%
}
\mapleinline{active}{1d}{\hspace{1.7cm}6*b*d-36*a*b+2*d+6*a*d=0,}{%
}
\mapleinline{active}{1d}{\hspace{1.7cm}-4*d^2+12*b*d-36*a*b+12*a*d=0,}{%
}
\mapleinline{active}{1d}{\hspace{1.7cm}-8*C-9+9*B+9*A-4*d+12*c=0,}{%
}
\mapleinline{active}{1d}{\hspace{1.7cm}-8*d-7+12*c+27*A+27*B-24*C=0,}{%
}
\mapleinline{active}{1d}{\hspace{1.7cm}8-3*a-3*b-32*C+27*A+27*B=0,}{%
}
\mapleinline{active}{1d}{\hspace{1.7cm}6-16*C+18*A+18*B-12*a+4*d-12*b=0,}{%
}
\mapleinline{active}{1d}{\hspace{1.7cm}4-12*a-12*b+8*d=0,}{%
}
\mapleinline{active}{1d}{\hspace{1.7cm}-C-2+3*c=0\};
}{%
}
\end{mapleinput}

\mapleresult
\begin{maplelatex}
\begin{eqnarray*}
\lefteqn{\mathit{LIST} := \{9\,B\,A + 4\,d - 6\,c\,d, \, - 9\,a\,
b + 9\,B\,A, } \\
 & &  - 18\,B\,A + 12\,d - 12\,c\,d + 4\,d^{2}, \,6\,b\,d - 36\,a
\,b + 2\,d + 6\,a\,d,  \\
 & &  - 4\,d^{2} + 12\,b\,d - 36\,a\,b + 12\,a\,d, \,4 - 12\,a -
12\,b + 8\,d,  \\
 & &  - C - 2 + 3\,c, \, - 8\,C - 9 + 9\,B + 9\,A - 4\,d + 12\,c
,  \\
 & &  - 8\,d - 7 + 12\,c + 27\,A + 27\,B - 24\,C,  \\
 & & 8 - 3\,a - 3\,b - 32\,C + 27\,A + 27\,B,  \\
 & & 6 - 16\,C + 18\,A + 18\,B - 12\,a + 4\,d - 12\,b\}
\;.
\end{eqnarray*}
\end{maplelatex}

\end{maplegroup}
\noindent
The {\tt solve} command gives the general solution:

\begin{maplegroup}
\begin{mapleinput}
\mapleinline{active}{1d}{solve(LIST,\{A,B,C,a,b,c,d\});}{%
}
\end{mapleinput}

\mapleresult
\begin{maplelatex}
\begin{eqnarray*}
\lefteqn{\{c={\displaystyle \frac {5}{6}}  + {\displaystyle 
\frac {1}{6}} \,d, \,A={\displaystyle \frac {1}{3}} \,d + 
{\displaystyle \frac {1}{3}} , \,B={\displaystyle \frac {1}{3}} 
\,d, \,b={\displaystyle \frac {1}{3}} \,d, \,a={\displaystyle 
\frac {1}{3}} \,d + {\displaystyle \frac {1}{3}} , \,C=
{\displaystyle \frac {1}{2}}  + {\displaystyle \frac {1}{2}} \,d
, \,d=d\}, } \\
 & & \{c={\displaystyle \frac {5}{6}}  + {\displaystyle \frac {1
}{6}} \,d, \,A={\displaystyle \frac {1}{3}} \,d + {\displaystyle 
\frac {1}{3}} , \,B={\displaystyle \frac {1}{3}} \,d, \,b=
{\displaystyle \frac {1}{3}} \,d + {\displaystyle \frac {1}{3}} 
, \,a={\displaystyle \frac {1}{3}} \,d, \,C={\displaystyle 
\frac {1}{2}}  + {\displaystyle \frac {1}{2}} \,d,  \\
 & & d=d\}, \{c={\displaystyle \frac {5}{6}}  + {\displaystyle 
\frac {1}{6}} \,d, \,b={\displaystyle \frac {1}{3}} \,d, \,B=
{\displaystyle \frac {1}{3}} \,d + {\displaystyle \frac {1}{3}} 
, \,A={\displaystyle \frac {1}{3}} \,d, \,a={\displaystyle 
\frac {1}{3}} \,d + {\displaystyle \frac {1}{3}} ,  \\
 & & C={\displaystyle \frac {1}{2}}  + {\displaystyle \frac {1}{2
}} \,d, \,d=d\}, \{c={\displaystyle \frac {5}{6}}  + 
{\displaystyle \frac {1}{6}} \,d, \,b={\displaystyle \frac {1}{3}
} \,d + {\displaystyle \frac {1}{3}} , \,B={\displaystyle \frac {
1}{3}} \,d + {\displaystyle \frac {1}{3}} , \,A={\displaystyle 
\frac {1}{3}} \,d,  \\
 & & a={\displaystyle \frac {1}{3}} \,d, \,C={\displaystyle 
\frac {1}{2}}  + {\displaystyle \frac {1}{2}} \,d, \,d=d\}
\;.
\mbox{\hspace{188pt}}
\end{eqnarray*}
\end{maplelatex}

\end{maplegroup}
\noindent
In an application, we will meet this example later again.

\subsection{Differentiation and Integration}

Differentiation is done using the differentiation rules. This is an
easy task. For our example function

\begin{maplegroup}
\begin{mapleinput}
\mapleinline{active}{1d}{input:=exp(x-x^2)*sin(x^6-1);}{%
}
\end{mapleinput}

\mapleresult
\begin{maplelatex}
\[
\mathit{input} := e^{(x - x^{2})}\,\mathrm{sin}(x^{6} - 1)
\]
\end{maplelatex}

\end{maplegroup}
\noindent
obviously the product rule is used:

\begin{maplegroup}
\begin{mapleinput}
\mapleinline{active}{1d}{derivative:=diff(input,x);}{%
}
\end{mapleinput}

\mapleresult
\begin{maplelatex}
\[
\mathit{derivative} := (1 - 2\,x)\,e^{(x - x^{2})}\,\mathrm{sin}(
x^{6} - 1) + 6\,e^{(x - x^{2})}\,\mathrm{cos}(x^{6} - 1)\,x^{5}
\;.
\]
\end{maplelatex}

\end{maplegroup}
\noindent
Integration is {\sl much} more difficult, and the different systems
have different approaches: Whereas Derive uses a good collection of
heuristics which enable the system to compute all explicitly given
integrals of Bronshtein and Semedyayev's integral table \cite{BS},
as already mentioned, Maple uses an algorithmic approach. 

In the sixties Risch developed an algorithm to compute an {\sl elementary}
antiderivative whenever one exists. If no such antiderivative exists,
his algorithm returns this information.
Here elementary means that both integrand and antiderivative are
rationally composed of exponentials and logarithms (see e.g.\ \cite{GCL}).
Adjoining the complex unit $i$ (denoted in Maple by $I$), 
trigonometric functions can be treated as well.

We integrate the derivative above. This takes a little longer:

\begin{maplegroup}
\begin{mapleinput}
\mapleinline{active}{1d}{integral:=int(derivative,x);}{%
}
\end{mapleinput}

\mapleresult
\begin{maplelatex}
\begin{eqnarray*}
\lefteqn{\mathit{integral} :=  - {\displaystyle \frac {1}{2}} \,I
\,e^{((x - 1)\,(I\,x^{5} + I\,x^{4} + I\,x^{3} + I\,x^{2} - x + I
\,x + I))}} \\
 & & \mbox{} + {\displaystyle \frac {1}{2}} \,I\,e^{( - (x - 1)\,
(I\,x^{5} + I\,x^{4} + I\,x^{3} + I\,x^{2} + x + I\,x + I))}
\;.
\mbox{\hspace{5pt}}
\end{eqnarray*}
\end{maplelatex}

\end{maplegroup}
\noindent
Since $i$ is adjoined, the resulting function looks not very familiar although
it is algebraically equal to our input function. In this particular
case, we can convert both functions to the same {\sl normal form} by
first converting exponentials to trigonometrics and applying then
rational factorization:

\begin{maplegroup}
\begin{mapleinput}
\mapleinline{active}{1d}{factor(convert(integral,trig));}{%
}
\end{mapleinput}

\mapleresult
\begin{maplelatex}
\begin{eqnarray*}
\lefteqn{ - \mathrm{sin}((x - 1)\,(x + 1)\,(x^{2} + x + 1)\,(x^{2
} - x + 1))} \\
 & & ( - \mathrm{cosh}(x\,(x - 1)) + \mathrm{sinh}(x\,(x - 1)))
\mbox{\hspace{7pt}}
\end{eqnarray*}
\end{maplelatex}

\end{maplegroup}
\begin{maplegroup}
\begin{mapleinput}
\mapleinline{active}{1d}{factor(convert(input,trig));}{%
}
\end{mapleinput}

\mapleresult
\begin{maplelatex}
\begin{eqnarray*}
\lefteqn{ - \mathrm{sin}((x - 1)\,(x + 1)\,(x^{2} + x + 1)\,(x^{2
} - x + 1))} \\
 & & ( - \mathrm{cosh}(x\,(x - 1)) + \mathrm{sinh}(x\,(x - 1)))
\mbox{\hspace{7pt}}
\end{eqnarray*}
\end{maplelatex}

\end{maplegroup}
\noindent
Note, however, that one can prove that for general transcendental
expressions a normal form does not exist.

\subsection{Differential Equations}

In engineering and in natural sciences the symbolic and numeric
solution of differential equations is rather important. We enter an
ordinary differential equation:

\begin{maplegroup}
\begin{mapleinput}
\mapleinline{active}{1d}{DE:=diff(y(x),x)=1+y(x)^2;}{%
}
\end{mapleinput}

\mapleresult
\begin{maplelatex}
\[
\mathit{DE} := {\frac {\partial }{\partial x}}\,\mathrm{y}(x)=1
 + \mathrm{y}(x)^{2}
\;.
\]
\end{maplelatex}

\end{maplegroup}
\noindent
After loading the {\tt DEtools} package, we can use the procedure
{\tt dfieldplot} to plot a direction field of the differential equation:

\begin{maplegroup}
\begin{mapleinput}
\mapleinline{active}{1d}{with(DEtools):}{%
}
\mapleinline{active}{1d}{dfieldplot(DE,y(x),x=-5..5,y=-5..5);}{%
}
\end{mapleinput}

\mapleresult
\begin{center}
\mapleplot{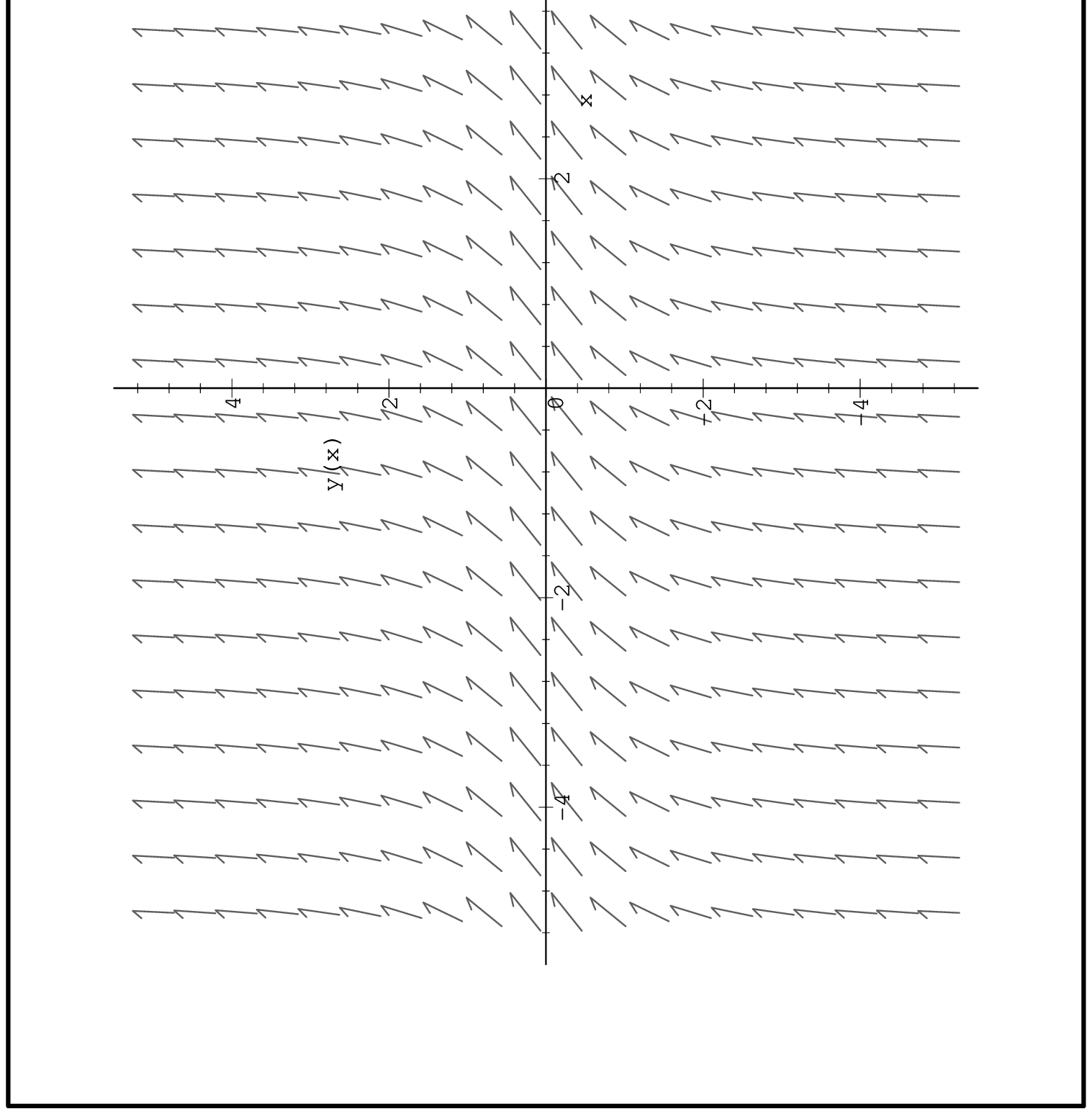}
\end{center}

\end{maplegroup}
\noindent
Given an initial value,
the command {\tt DEplot} plots a numeric solution by a Runge-Kutta 
type approach:

\begin{maplegroup}
\begin{mapleinput}
\mapleinline{active}{1d}{DEplot(\{DE\},\{y(x)\},x=-1..1,[[y(0)=0]]);}{
}
\end{mapleinput}

\mapleresult
\begin{center}
\mapleplot{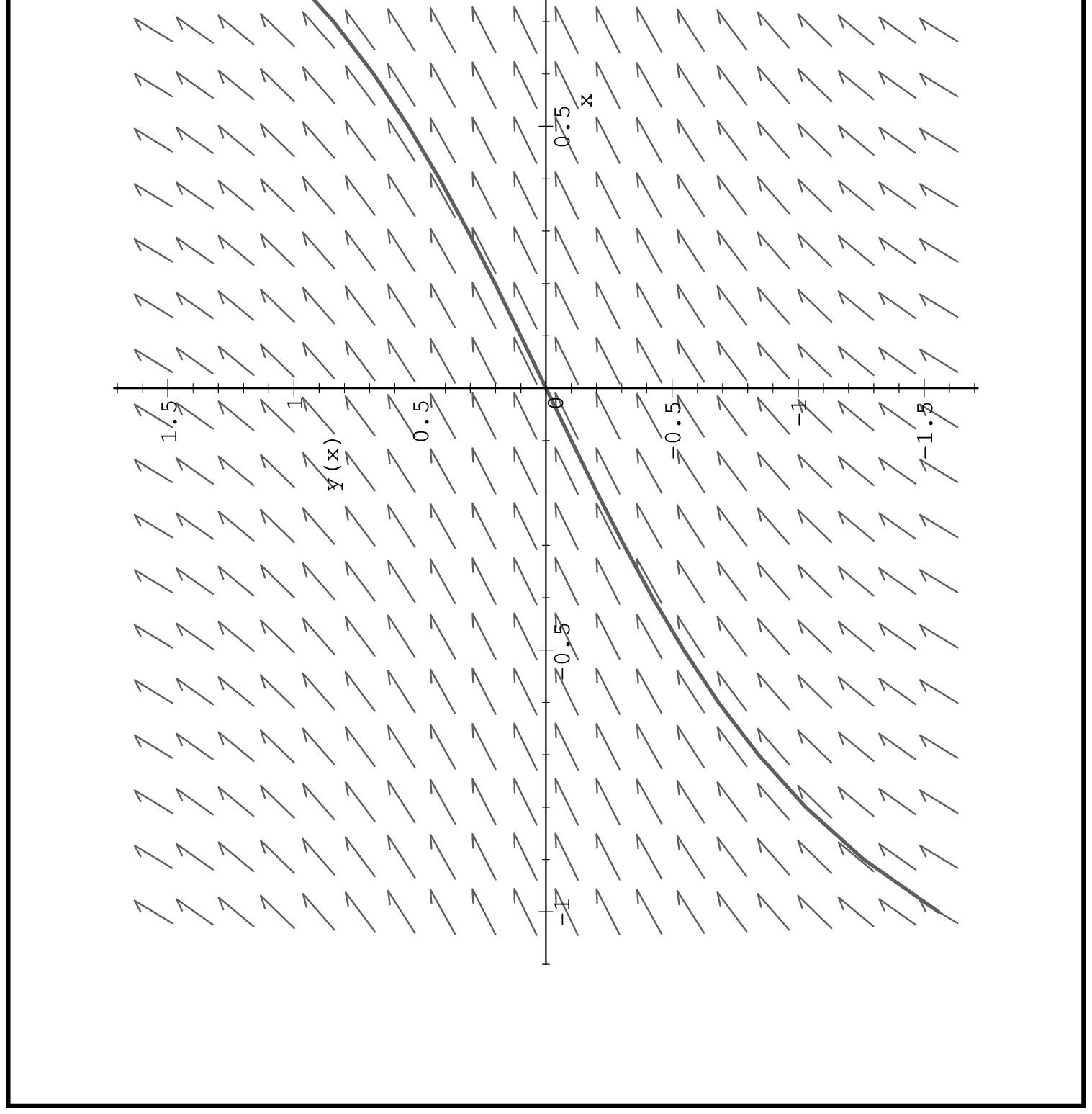}
\end{center}

\end{maplegroup}
\noindent
Using a combination of heuristic and algorithmic techniques, Maple can
solve many ordinary differential equations explicitly. Our 
initial value problem has the solution:

\begin{maplegroup}
\begin{mapleinput}
\mapleinline{active}{1d}{dsolve(\{DE,y(0)=0\},y(x));}{%
}
\end{mapleinput}

\mapleresult
\begin{maplelatex}
\[
\mathrm{y}(x)=\mathrm{tan}(x)
\;.
\]
\end{maplelatex}

\end{maplegroup}
\noindent
As another example, we consider a linear differential equation of second order.

\begin{maplegroup}
\begin{mapleinput}
\mapleinline{active}{1d}{DE:=diff(y(x),x$2)-y(x)=sin(x)*x;}{%
}
\end{mapleinput}

\mapleresult
\begin{maplelatex}
\[
\mathit{DE} := ({\frac {\partial ^{2}}{\partial x^{2}}}\,\mathrm{
y}(x)) - \mathrm{y}(x)=\mathrm{sin}(x)\,x
\]
\end{maplelatex}

\end{maplegroup}
\noindent
with explicit solution

\begin{maplegroup}
\begin{mapleinput}
\mapleinline{active}{1d}{dsolve(DE,y(x));}{%
}
\end{mapleinput}

\mapleresult
\begin{maplelatex}
\begin{eqnarray*}
\lefteqn{\mathrm{y}(x)=({\displaystyle \frac {1}{2}} \,( - 
{\displaystyle \frac {1}{2}} \,x + {\displaystyle \frac {1}{2}} )
\,e^{x}\,\mathrm{cos}(x) + {\displaystyle \frac {1}{4}} \,
\mathrm{sin}(x)\,x\,e^{x} + {\displaystyle \frac {1}{2}} \,( - 
{\displaystyle \frac {1}{2}} \,x - {\displaystyle \frac {1}{2}} )
\,e^{( - x)}\,\mathrm{cos}(x)} \\
 & & \mbox{} - {\displaystyle \frac {1}{4}} \,x\,e^{( - x)}\,
\mathrm{sin}(x))\mathrm{sinh}(x)\mbox{} + ( - {\displaystyle 
\frac {1}{2}} \,( - {\displaystyle \frac {1}{2}} \,x + 
{\displaystyle \frac {1}{2}} )\,e^{x}\,\mathrm{cos}(x) - 
{\displaystyle \frac {1}{4}} \,\mathrm{sin}(x)\,x\,e^{x} \\
 & & \mbox{} + {\displaystyle \frac {1}{2}} \,( - {\displaystyle 
\frac {1}{2}} \,x - {\displaystyle \frac {1}{2}} )\,e^{( - x)}\,
\mathrm{cos}(x) - {\displaystyle \frac {1}{4}} \,x\,e^{( - x)}\,
\mathrm{sin}(x))\mathrm{cosh}(x)\mbox{} + \mathit{\_C1}\,\mathrm{
sinh}(x) \\
 & & \mbox{} + \mathit{\_C2}\,\mathrm{cosh}(x)
\;.
\end{eqnarray*}
\end{maplelatex}

\end{maplegroup}
\noindent
A plot based on a numerical computation is given by

\begin{maplegroup}
\begin{mapleinput}
\mapleinline{active}{1d}{DEplot(\{DE\},\{y(x)\},x=-5..5,[[y(0)=0,D(y)(0)=1]]);}{%
}
\end{mapleinput}

\mapleresult
\begin{center}
\mapleplot{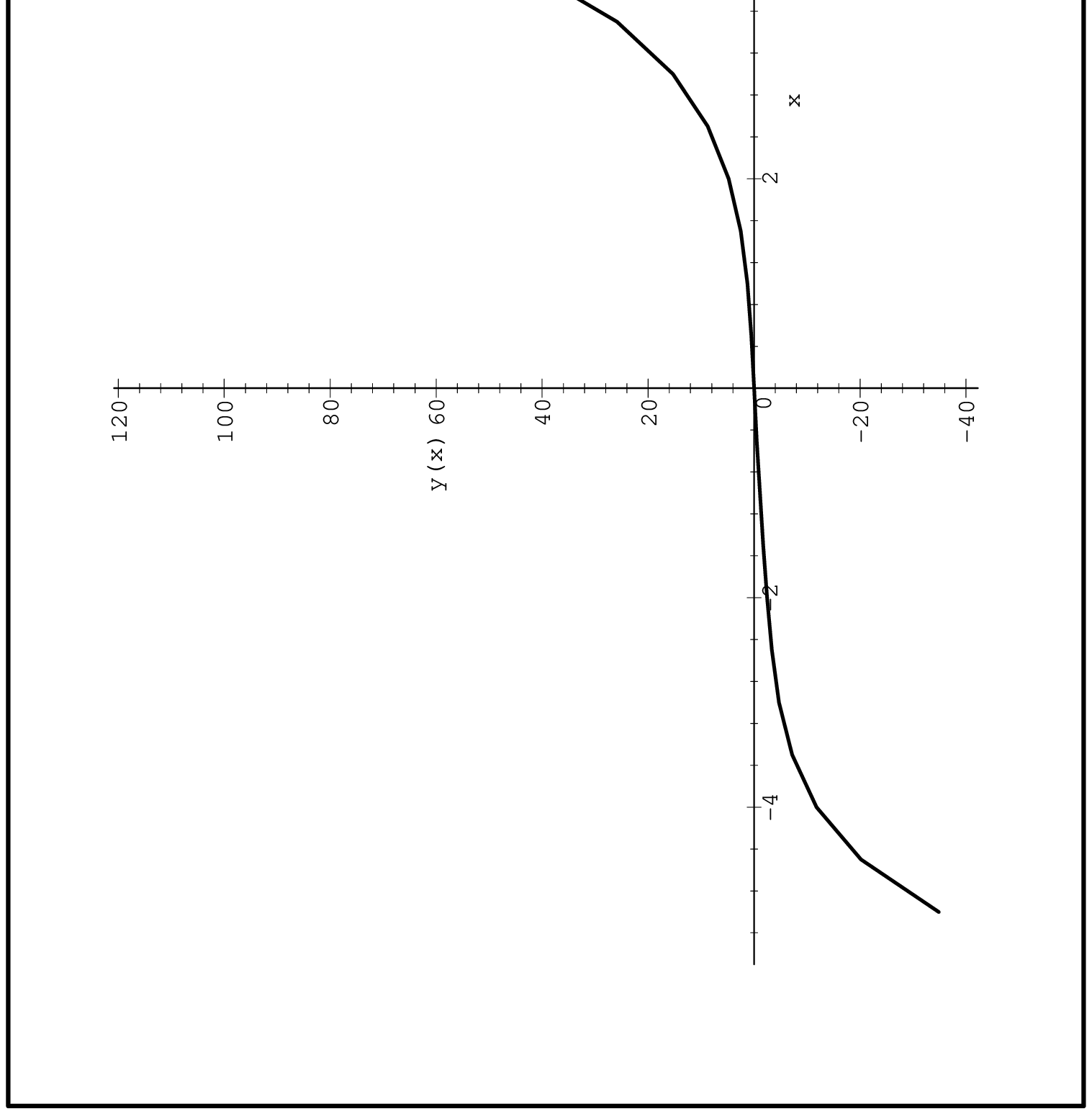}
\end{center}

\end{maplegroup}
\noindent
The corresponding initial value problem has the explicit solution

\begin{maplegroup}
\begin{mapleinput}
\mapleinline{active}{1d}{solution:=dsolve(\{DE,y(0)=0,D(y)(0)=1\},y(x));}{%
}
\end{mapleinput}

\mapleresult
\begin{maplelatex}
\begin{eqnarray*}
\lefteqn{\mathit{solution} := \mathrm{y}(x)=({\displaystyle 
\frac {1}{2}} \,( - {\displaystyle \frac {1}{2}} \,x + 
{\displaystyle \frac {1}{2}} )\,e^{x}\,\mathrm{cos}(x) + 
{\displaystyle \frac {1}{4}} \,\mathrm{sin}(x)\,x\,e^{x}} \\
 & & \mbox{} + {\displaystyle \frac {1}{2}} \,( - {\displaystyle 
\frac {1}{2}} \,x - {\displaystyle \frac {1}{2}} )\,e^{( - x)}\,
\mathrm{cos}(x) - {\displaystyle \frac {1}{4}} \,x\,e^{( - x)}\,
\mathrm{sin}(x))\mathrm{sinh}(x)\mbox{} + ( \\
 & &  - {\displaystyle \frac {1}{2}} \,( - {\displaystyle \frac {
1}{2}} \,x + {\displaystyle \frac {1}{2}} )\,e^{x}\,\mathrm{cos}(
x) - {\displaystyle \frac {1}{4}} \,\mathrm{sin}(x)\,x\,e^{x} + 
{\displaystyle \frac {1}{2}} \,( - {\displaystyle \frac {1}{2}} 
\,x - {\displaystyle \frac {1}{2}} )\,e^{( - x)}\,\mathrm{cos}(x)
 \\
 & & \mbox{} - {\displaystyle \frac {1}{4}} \,x\,e^{( - x)}\,
\mathrm{sin}(x))\mathrm{cosh}(x)\mbox{} + \mathrm{sinh}(x) + 
{\displaystyle \frac {1}{2}} \,\mathrm{cosh}(x)
\end{eqnarray*}
\end{maplelatex}

\end{maplegroup}
\noindent
which can be simplified to

\begin{maplegroup}
\begin{mapleinput}
\mapleinline{active}{1d}{simplify(convert(rhs(solution),trig));}{%
}
\end{mapleinput}

\mapleresult
\begin{maplelatex}
\[
\mathrm{sinh}(x) + {\displaystyle \frac {1}{2}} \,\mathrm{cosh}(x
) - {\displaystyle \frac {1}{2}} \,\mathrm{cos}(x) - 
{\displaystyle \frac {1}{2}} \,\mathrm{sin}(x)\,x
\;.
\]
\end{maplelatex}

\end{maplegroup}

\subsection{Formal Power Series and Differential Equations}
Next, we consider the opposite problem to generate differential equations from
expressions. This will lead us also to the generation of power series
of hypergeometric type.

After loading the {\tt FPS} package \cite{GK}

\begin{maplegroup}
\begin{mapleinput}
\mapleinline{active}{1d}{with(share): with(FPS):}{%
}
\end{mapleinput}

\mapleresult
\begin{maplettyout}
See ?share and ?share,contents for information about the share library
Share Library:  FPS
Author: Gruntz, Dominik.
Description:  FPS function attempts to find a formal power
series expansion for a function in terms of a formula for the
coefficients
\end{maplettyout}

\end{maplegroup}
\noindent
we can, e.g., compute the formal power series of the square of the inverse 
tangent function:

\begin{maplegroup}
\begin{mapleinput}
\mapleinline{active}{1d}{FPS(arcsin(x)^2,x);}{%
}
\end{mapleinput}

\mapleresult
\begin{maplelatex}
\begin{equation}
{\displaystyle \sum _{k=0}^{\infty }} \,{\displaystyle \frac {(k
\mathrm{!})^{2}\,4^{k}\,x^{(2\,k + 2)}}{(k + 1)\,(1 + 2\,k)
\mathrm{!}}} 
\;.
\label{eq:FPS}
\end{equation}
\end{maplelatex}

\end{maplegroup}
\noindent
The algorithm behind this procedure is the following (\cite{Koepf1990},
\cite{GK}):

In the first step, by linear algebra techniques,
a homogeneous linear differential equation with polynomial coefficients
is sought for the given expression

\begin{maplegroup}
\begin{mapleinput}
\mapleinline{active}{1d}{DE:=SimpleDE(arcsin(x)^2,x,F);}{%
}
\end{mapleinput}

\mapleresult
\begin{maplelatex}
\[
\mathit{DE} := (x - 1)\,(x + 1)\,({\frac {\partial ^{3}}{
\partial x^{3}}}\,\mathrm{F}(x)) + ({\frac {\partial }{\partial x
}}\,\mathrm{F}(x)) + 3\,x\,({\frac {\partial ^{2}}{\partial x^{2}
}}\,\mathrm{F}(x))=0
\;.
\]
\end{maplelatex}

\end{maplegroup}
\noindent
We call such a differential equation as well as the corresponding function
{\sl holonomic}. Next, substituting
the series 
\[
F(x)=\sum_{k=0}^\infty a_k\, x^k
\]
in this differential equation and equating coefficients yields
the holonomic recurrence equation for $a_k$:

\begin{maplegroup}
\begin{mapleinput}
\mapleinline{active}{1d}{RE:=SimpleRE(arcsin(x)^2,x,a);}{%
}
\end{mapleinput}

\mapleresult
\begin{maplelatex}
\begin{eqnarray*}
\lefteqn{\mathit{RE} :=  - (k + 1)(\mathrm{a}(k + 3)\,k^{2} - k^{
2}\,\mathrm{a}(k + 1) + 5\,\mathrm{a}(k + 3)\,k - 2\,k\,\mathrm{a
}(k + 1)} \\
 & & \mbox{} - \mathrm{a}(k + 1) + 6\,\mathrm{a}(k + 3))=0
\mbox{\hspace{154pt}}
\end{eqnarray*}
\end{maplelatex}

\end{maplegroup}
\noindent
which can be put in factored form

\begin{maplegroup}
\begin{mapleinput}
\mapleinline{active}{1d}{map(factor,collect(lhs(RE),a))=0;}{%
}
\end{mapleinput}

\mapleresult
\begin{maplelatex}
\begin{equation}
- (k + 1)\,(k + 2)\,(k + 3)\,\mathrm{a}(k + 3) + (k + 1)^{3}\,
\mathrm{a}(k + 1)=0
\;.
\label{eq:ak}
\end{equation}
\end{maplelatex}

\end{maplegroup}
\noindent
Notice that the resulting recurrence equation gives $a_{k+2}$ as a
rational multiple of $a_k$. If $A_{k+1}$ is a
rational multiple of $A_k$ then it is called a {\sl hypergeometric term}.
From (\ref{eq:ak}), $a_k$ can be easily computed using two initial values.
This finally generates
the explicit series representation (\ref{eq:FPS}). Note, however, that for
an explicit representation the above factorization is necessary;
see (\ref{eq:pFq}).

\noindent

By solving the differential equation for $F(x)=\arcsin^2(x)$ with two initial
values, we would like to reconstruct the input:

\begin{maplegroup}

\begin{mapleinput}
\mapleinline{active}{1d}{solution:=dsolve(\{DE,F(0)=0,D(F)(0)=0,(D@@2)(F)(0)=2\},F(x));}{%
}
\end{mapleinput}

\mapleresult
\begin{maplelatex}
\begin{eqnarray*}
\lefteqn{\mathit{solution} := } \\
 & & \mathrm{F}(x)={\displaystyle \frac {1}{4}} \,\pi ^{2} + I\,
\pi \,\mathrm{ln}(x + \sqrt{(x - 1)\,(x + 1)}) - \mathrm{ln}(x + 
\sqrt{(x - 1)\,(x + 1)})^{2}
\end{eqnarray*}
\end{maplelatex}

\end{maplegroup}
\begin{maplegroup}
\begin{mapleinput}
\mapleinline{active}{1d}{convert(arcsin(x)^2,ln);}{%
}
\end{mapleinput}

\mapleresult
\begin{maplelatex}
\[
 - \mathrm{ln}(\sqrt{1 - x^{2}} + I\,x)^{2}
\;.
\]
\end{maplelatex}

\end{maplegroup}
\noindent
As before, we see that transcendental functions come in quite 
different disguises.

It turns out that sum and product of two holonomic functions are again
holonomic, and the corresponding holonomic (differential or recurrence)
equations can be constructed from the given holonomic equations by linear 
algebra (\cite{Beke1}, \cite{Beke2}, \cite{Sta}, \cite{SZ}).

As an example, we consider both the sum and the product of the functions 
$f(x)=\arcsin x$ and $g(x)=e^x$. Here are their holonomic equations:

\begin{maplegroup}

\begin{mapleinput}
\mapleinline{active}{1d}{DE1:=SimpleDE(arcsin(x),x,F);}{%
}
\end{mapleinput}

\mapleresult
\begin{maplelatex}
\[
\mathit{DE1} := (x - 1)\,(x + 1)\,({\frac {\partial ^{2}}{
\partial x^{2}}}\,\mathrm{F}(x)) + ({\frac {\partial }{\partial x
}}\,\mathrm{F}(x))\,x=0
\]
\end{maplelatex}

\end{maplegroup}
\begin{maplegroup}
\begin{mapleinput}
\mapleinline{active}{1d}{DE2:=SimpleDE(exp(x),x,F);}{%
}
\end{mapleinput}

\mapleresult
\begin{maplelatex}
\[
\mathit{DE2} := ({\frac {\partial }{\partial x}}\,\mathrm{F}(x))
 - \mathrm{F}(x)=0
\;.
\]
\end{maplelatex}

\end{maplegroup}
\noindent
From these, we can compute the holonomic equations that are valid
for $f(x)+g(x)$ and $f(x)\cdot g(x)$. For this purpose, we load the 
{\tt gfun} package \cite{SZ}:

\begin{maplegroup}
\begin{mapleinput}
\mapleinline{active}{1d}{with(gfun);}{%
}
\end{mapleinput}

\mapleresult
\begin{maplelatex}
\begin{eqnarray*}
\lefteqn{[\mathit{Laplace}, \,\mathit{algebraicsubs}, \,\mathit{
algeqtodiffeq}, \,\mathit{algeqtoseries}, \,\mathit{algfuntoalgeq
}, \,\mathit{borel}, } \\
 & & \mathit{cauchyproduct}, \,\mathit{diffeq*diffeq}, \,\mathit{
diffeq+diffeq}, \,\mathit{diffeqtorec}, \,\mathit{guesseqn},  \\
 & & \mathit{guessgf}, \,\mathit{hadamardproduct}, \,\mathit{
holexprtodiffeq}, \,\mathit{invborel}, \,\mathit{listtoalgeq}, 
 \\
 & & \mathit{listtodiffeq}, \,\mathit{listtohypergeom}, \,
\mathit{listtolist}, \,\mathit{listtoratpoly}, \,\mathit{
listtorec},  \\
 & & \mathit{listtoseries}, \,\mathit{listtoseries/Laplace}, \,
\mathit{listtoseries/egf}, \,\mathit{listtoseries/lgdegf},  \\
 & & \mathit{listtoseries/lgdogf}, \,\mathit{listtoseries/ogf}, 
\,\mathit{listtoseries/revegf},  \\
 & & \mathit{listtoseries/revogf}, \,\mathit{maxdegcoeff}, \,
\mathit{maxdegeqn}, \,\mathit{maxordereqn},  \\
 & & \mathit{mindegcoeff}, \,\mathit{mindegeqn}, \,\mathit{
minordereqn}, \,\mathit{optionsgf}, \,\mathit{poltodiffeq},  \\
 & & \mathit{poltorec}, \,\mathit{ratpolytocoeff}, \,\mathit{
rec*rec}, \,\mathit{rec+rec}, \,\mathit{rectodiffeq}, \,\mathit{
rectoproc},  \\
 & & \mathit{seriestoalgeq}, \,\mathit{seriestodiffeq}, \,
\mathit{seriestohypergeom}, \,\mathit{seriestolist},  \\
 & & \mathit{seriestoratpoly}, \,\mathit{seriestorec}, \,\mathit{
seriestoseries}]\mbox{\hspace{118pt}}
\end{eqnarray*}
\end{maplelatex}

\end{maplegroup}
\noindent
The procedures \verb@`diffeq+diffeq`@ and \verb@`diffeq*diffeq`@
compute the differential equations of sum and product, respectively:

\begin{maplegroup}
\begin{mapleinput}
\mapleinline{active}{1d}{`diffeq+diffeq`(DE1,DE2,F(x));}{%
}
\end{mapleinput}

\mapleresult
\begin{maplelatex}
\begin{eqnarray*}
\lefteqn{\{( - x^{3} - 2\,x^{2} + x - 1)\,\mathrm{D}(F)(x) + ( - 
x^{4} + 4\,x^{2})\,(\mathrm{D}^{(2)})(F)(x)} \\
 & & \mbox{} + (1 - 2\,x^{2} + x^{4} - x + x^{3})\,(\mathrm{D}^{(
3)})(F)(x), \,(\mathrm{D}^{(2)})(F)(0)={\mathit{\_C}_{0}}\}
\end{eqnarray*}
\end{maplelatex}

\end{maplegroup}
\begin{maplegroup}
\begin{mapleinput}
\mapleinline{active}{1d}{`diffeq*diffeq`(DE1,DE2,F(x));}{%
}
\end{mapleinput}

\mapleresult
\begin{maplelatex}
\[
( - 1 + x^{2} - x)\,\mathrm{F}(x) + (x + 2 - 2\,x^{2})\,\mathrm{D
}(F)(x) + ( - 1 + x^{2})\,(\mathrm{D}^{(2)})(F)(x)
\]
\end{maplelatex}

\end{maplegroup}
\noindent
which we could also have obtained using {\tt SimpleDE} directly:%
\footnote{Note that {\tt SimpleDE} uses a slightly different approach 
(also based on linear algebra) that
sometimes can find differential equations of lower order than
{\tt `diffeq+diffeq`} and {\tt `diffeq*diffeq`}.}

\begin{maplegroup}

\begin{mapleinput}
\mapleinline{active}{1d}{SimpleDE(arcsin(x)+exp(x),x,F);}{%
}
\end{mapleinput}

\mapleresult
\begin{maplelatex}
\begin{eqnarray*}
\lefteqn{(x - 1)\,(x + 1)\,(x - 1 + x^{2})\,({\frac {\partial ^{3
}}{\partial x^{3}}}\,\mathrm{F}(x)) + (x - x^{3} - 2\,x^{2} - 1)
\,({\frac {\partial }{\partial x}}\,\mathrm{F}(x))} \\
 & & \mbox{} - x^{2}\,(x - 2)\,(x + 2)\,({\frac {\partial ^{2}}{
\partial x^{2}}}\,\mathrm{F}(x))=0\mbox{\hspace{119pt}}
\end{eqnarray*}
\end{maplelatex}

\end{maplegroup}
\begin{maplegroup}
\begin{mapleinput}
\mapleinline{active}{1d}{SimpleDE(arcsin(x)*exp(x),x,F);}{%
}
\end{mapleinput}

\mapleresult
\begin{maplelatex}
\begin{eqnarray*}
\lefteqn{(x - 1)\,(x + 1)\,({\frac {\partial ^{2}
}{\partial x^{2}}}\,\mathrm{F}(x)) + (x + 2 - 2\,x^{2})\,(
{\frac {\partial }{\partial x}}\,\mathrm{F}(x))} \\
 & & \mbox{} + ( - 1 + x^{2} - x)\,\mathrm{F}(x)=0
\;.
\mbox{\hspace{118pt}}
\end{eqnarray*}
\end{maplelatex}

\end{maplegroup}

\noindent
{\tt SimpleDE}
can also generate differential equations for some special functions,
e.g., for the Bessel functions $J_n(x)$:

\begin{maplegroup}

\begin{mapleinput}
\mapleinline{active}{1d}{DE:=SimpleDE(BesselJ(n,x),x,F);}{%
}
\end{mapleinput}

\mapleresult
\begin{maplelatex}
\[
\mathit{DE} := ({\frac {\partial ^{2}}{\partial x^{2}}}\,\mathrm{
F}(x))\,x^{2} - (n - x)\,(n + x)\,\mathrm{F}(x) + ({\frac {
\partial }{\partial x}}\,\mathrm{F}(x))\,x=0
\;.
\]
\end{maplelatex}

\end{maplegroup}
\noindent
Maple can solve this differential equation easily:

\begin{maplegroup}
\begin{mapleinput}
\mapleinline{active}{1d}{dsolve(DE,F(x));}{%
}
\end{mapleinput}

\mapleresult
\begin{maplelatex}
\[
\mathrm{F}(x)=\mathit{\_C1}\,\mathrm{BesselY}(n, \,x) + \mathit{
\_C2}\,\mathrm{BesselJ}(n, \,x)
\;.
\]
\end{maplelatex}

\end{maplegroup}
\noindent
Even the more complicated differential equation of the product

\begin{maplegroup}
\begin{mapleinput}
\mapleinline{active}{1d}{DE:=SimpleDE(BesselJ(n,x)*exp(x),x,F);}{%
}
\end{mapleinput}

\mapleresult
\begin{maplelatex}
\[
\mathit{DE} := (2\,x^{2} - n^{2} - x)\,\mathrm{F}(x) - ( - 1 + 2
\,x)\,x\,({\frac {\partial }{\partial x}}\,\mathrm{F}(x)) + (
{\frac {\partial ^{2}}{\partial x^{2}}}\,\mathrm{F}(x))\,x^{2}=0
\]
\end{maplelatex}

\end{maplegroup}
\noindent
can be treated by Maple

\begin{maplegroup}
\begin{mapleinput}
\mapleinline{active}{1d}{dsolve(DE,F(x));}{%
}
\end{mapleinput}

\mapleresult
\begin{maplelatex}
\[
\mathrm{F}(x)=\mathit{\_C1}\,\mathrm{BesselJ}(n, \,x)\,e^{x} + 
\mathit{\_C2}\,\mathrm{BesselY}(n, \,x)\,e^{x}
\;,
\]
\end{maplelatex}

\end{maplegroup}
\noindent
but for the differential equation

\begin{maplegroup}
\begin{mapleinput}
\mapleinline{active}{1d}{DE:=SimpleDE(BesselJ(n,x)+exp(x),x,F);}{%
}
\end{mapleinput}

\mapleresult
\begin{maplelatex}
\begin{eqnarray*}
\lefteqn{\mathit{DE} := (2\,x^{4} - 3\,x^{2}\,n^{2} + x^{3} - 3\,
n^{2}\,x - x^{2} + n^{4} - n^{2})\,\mathrm{F}(x)} \\
 & & \mbox{} + ( - n^{4} + 3\,x^{2} + n^{2} + x^{3} - 2\,x^{4} + 
3\,x^{2}\,n^{2})\,({\frac {\partial }{\partial x}}\,\mathrm{F}(x)
) \\
 & & \mbox{} - ( - 2\,x^{3} + n^{2}\,x + x^{2} - 3\,n^{2} + 2\,x)
\,x\,({\frac {\partial ^{2}}{\partial x^{2}}}\,\mathrm{F}(x)) \\
 & & \mbox{} + x^{2}\,( - 2\,x^{2} + n^{2} - x)\,({\frac {
\partial ^{3}}{\partial x^{3}}}\,\mathrm{F}(x))=0
\;,
\end{eqnarray*}
\end{maplelatex}

\end{maplegroup}
\noindent
Maple fails:

\begin{maplegroup}
\begin{mapleinput}
\mapleinline{active}{1d}{dsolve(DE,F(x));}{%
}
\end{mapleinput}

\mapleresult
\begin{maplelatex}
\begin{eqnarray*}
\lefteqn{\mathrm{F}(x)=\mathit{\_C1}\,e^{x} + e^{x}\mathrm{DESol}
 \left( {\vrule height1.38em width0em depth1.38em} \right. \! 
 \! } \\
 & &  \left\{  \! \mathrm{\_Y}(x) - {\displaystyle \frac {x\,(2\,
x + 1)\,({\frac {\partial }{\partial x}}\,\mathrm{\_Y}(x))}{ - 2
\,x^{2} + n^{2} - x}}  - {\displaystyle \frac {( - 2\,x^{4} + x^{
2}\,n^{2} - x^{3})\,({\frac {\partial ^{2}}{\partial x^{2}}}\,
\mathrm{\_Y}(x))}{( - 2\,x^{2} + n^{2} - x)^{2}}}  \!  \right\} 
,  \\
 & & \{\mathrm{\_Y}(x)\} \! \! \left. {\vrule 
height1.38em width0em depth1.38em} \right) 
\end{eqnarray*}
\end{maplelatex}

\end{maplegroup}
\begin{maplegroup}
\begin{mapleinput}
\end{mapleinput}

\end{maplegroup}
\noindent
although Maple was able to find the exponential summand (and hence
reduced the order by one). 

Nevertheless,
it is not astonishing that Maple cannot find all such solutions since
for this type of nonelementary solutions no algorithms exist.

\section{Special Functions and Computer Algebra}

Power series of hypergeometric type---the example function $\arcsin^2 x$ as
well as the Bessel functions are of this type, e.g.---are the most important 
special functions.

The {\sl generalized hypergeometric series} is given by
\begin{equation}
_{p}F_{q}\left.\left(\begin{array}{cccc}
a_{1}&a_{2}&\cdots&a_{p}\\
b_{1}&b_{2}&\cdots&b_{q}\\
            \end{array}\right| x\right)
:=
\sum_{k=0}^\infty A_k\,x^{k}
=\sum_{k=0}^\infty \frac
{(a_{1})_{k}\cdot(a_{2})_{k}\cdots(a_{p})_{k}}
{(b_{1})_{k}\cdot(b_{2})_{k}\cdots(b_{q})_{k}\,k!}x^{k}
\label{eq:pFq}
\end{equation}
where $(a)_{k}:=\prod\limits_{j=1}^k (a\!+\!j\!-\!1)=\Gamma (a+k)/\Gamma(a)$
denotes the {\sl Pochhammer-Symbol} or {\sl shifted factorial}.

$A_k$ is a hypergeometric term and fulfils the
recurrence equation $(k\in\N)$
\[
A_{k+1}:=\frac{(k+a_1)\cdots (k+a_p)}{(k+b_1)\cdots
(k+b_q)(k+1)}\cdot A_k
\]
with the initial value
\[
A_0:=1\;.
\]
In Maple the hypergeometric series is given as
{\tt hypergeom(plist,qlist,x)}, where
\[
{\rm plist}=[a_{1},a_{2},\ldots,a_p]\;,
\]
\[
{\rm qlist}=[b_{1},b_{2},\ldots,b_q]\;.
\]
Here are some more hypergeometric examples:

\begin{maplegroup}

\begin{mapleinput}
\mapleinline{active}{1d}{F:=sqrt(x)*arcsin(sqrt(x))+sqrt(1-x);}{%
}
\end{mapleinput}

\mapleresult
\begin{maplelatex}
\[
F := \sqrt{x}\,\mathrm{arcsin}(\sqrt{x}) + \sqrt{1 - x}
\]
\end{maplelatex}

\end{maplegroup}
\begin{maplegroup}
\begin{mapleinput}
\mapleinline{active}{1d}{SUM:=FPS(F,x);}{%
}
\end{mapleinput}

\mapleresult
\begin{maplelatex}
\[
\mathit{SUM} := {\displaystyle \sum _{k=0}^{\infty }} \,
{\displaystyle \frac {4^{( - k)}\,(2\,k)\mathrm{!}\,x^{k}}{(k
\mathrm{!})^{2}\,(2\,k - 1)^{2}}} 
\]
\end{maplelatex}

\end{maplegroup}
\begin{maplegroup}
\begin{mapleinput}
\mapleinline{active}{1d}{convert(SUM,hypergeom);}{%
}
\end{mapleinput}

\mapleresult
\begin{maplelatex}
\[
\mathrm{hypergeom}([{\displaystyle \frac {-1}{2}} , \,
{\displaystyle \frac {-1}{2}} ], \,[{\displaystyle \frac {1}{2}} 
], \,x)
\]
\end{maplelatex}

\end{maplegroup}
\begin{maplegroup}
\begin{mapleinput}

\mapleinline{active}{1d}{F:=-(sqrt(Pi)/2*sqrt(x)*erf(sqrt(x))*(1+1/2/x)+exp(-x)/2);
}{%
}
\end{mapleinput}

\mapleresult
\begin{maplelatex}
\[
F :=  - {\displaystyle \frac {1}{2}} \,\sqrt{\pi }\,\sqrt{x}\,
\mathrm{erf}(\sqrt{x})\,(1 + {\displaystyle \frac {1}{2}} \,
{\displaystyle \frac {1}{x}} ) - {\displaystyle \frac {1}{2}} \,e
^{( - x)}
\]
\end{maplelatex}

\end{maplegroup}
\begin{maplegroup}
\begin{mapleinput}
\mapleinline{active}{1d}{SUM:=FPS(F,x);}{%
}
\end{mapleinput}

\mapleresult
\begin{maplelatex}
\[
\mathit{SUM} := {\displaystyle \sum _{k=0}^{\infty }} \,
{\displaystyle \frac {(-1)^{k}\,x^{k}}{k\mathrm{!}\,(2\,k + 1)\,(
 - 1 + 2\,k)}} 
\]
\end{maplelatex}

\end{maplegroup}
\begin{maplegroup}
\begin{mapleinput}
\mapleinline{active}{1d}{convert(SUM,hypergeom);}{%
}
\end{mapleinput}

\mapleresult
\begin{maplelatex}
\[
 - \mathrm{KummerM}({\displaystyle \frac {-1}{2}} , \,
{\displaystyle \frac {3}{2}} , \, - x)
\]
\end{maplelatex}

\end{maplegroup}
\noindent
With {\tt convert}, one can convert series into hypergeometric notation;
{\tt KummerM} is another name for the confluent hypergeometric
function $_1 F_1$.

\subsection{Summation}

Whereas the {\tt FPS} command converts expressions into series representations,
the opposite question is to find explicit representations for sums. Note
that the examples of the remaining paper are from \cite{Koepf}.

The main interest lies in sums of hypergeometric terms. 
As an example, we ask: Why does Maple evaluate the sum

\begin{maplegroup}

\begin{mapleinput}
\mapleinline{active}{1d}{sum((-1)^k*binomial(n,k),k=a..b);}{%
}
\end{mapleinput}

\mapleresult
\begin{maplelatex}
\[
 - {\displaystyle \frac {(b + 1)\,(-1)^{(b + 1)}\,\mathrm{
binomial}(n, \,b + 1)}{n}}  + {\displaystyle \frac {a\,(-1)^{a}\,
\mathrm{binomial}(n, \,a)}{n}} 
\]
\end{maplelatex}

\end{maplegroup}
\noindent
for arbitrary bounds $a$ and $b$ in simple form, but fails with

\begin{maplegroup}
\begin{mapleinput}
\mapleinline{active}{1d}{sum(binomial(n,k),k=a..b);}{%
}
\end{mapleinput}

\mapleresult
\begin{maplelatex}
\begin{eqnarray*}
\lefteqn{\mathrm{binomial}(n, \,a)\,\mathrm{hypergeom}([1, \, - n
 + a], \,[1 + a], \,-1)} \\
 & & \mbox{} - \mathrm{binomial}(n, \,b + 1)\,\mathrm{hypergeom}(
[1, \, - n + b + 1], \,[2 + b], \,-1)
\;?
\end{eqnarray*}
\end{maplelatex}

\end{maplegroup}
\noindent
On the other hand, for the special bounds $a=0$ and $b=n$, Maple is successful,
again:

\begin{maplegroup}
\begin{mapleinput}
\mapleinline{active}{1d}{sum(binomial(n,k),k=0..n);}{%
}
\end{mapleinput}

\mapleresult
\begin{maplelatex}
\[
2^{n}
\;.
\]
\end{maplelatex}

\end{maplegroup}
\noindent
The reason for this behavior is that the first summand, $(-1)^k\,{n\choose k}$,
has a hypergeometric term antidifference (w.r.t.\ the variable $k$), 
and the second one, ${n\choose k}$,
has not. The last sum is a definite sum with natural bounds,
i.e., the sum can be considered as infinite sum $(k=-\infty..\infty)$,
and the result, again, is a hypergeometric term (w.r.t.\ the variable $n$).
We will see how we can find these types of results algorithmically.

\pagebreak

$s_k$ is called an {\sl antidifference} of $a_k$, if
\[
s_{k+1}-s_k=a_k
\;,
\]
If such an antidifference is known, then summation is trivial since
by telescoping
\[
\sum_{k=a}^b a_k=(s_{b+1}-s_b)+(s_{b}-s_{b-1})+\cdots+(s_{a+1}-s_a)
=s_{b+1}-s_a
\;.
\]
This is very similar to the integration case.

The antidifference of the first summand is given by
\begin{maplegroup}
\begin{mapleinput}
\mapleinline{active}{1d}{sum((-1)^k*binomial(n,k),k);}{%
}
\end{mapleinput}

\mapleresult
\begin{maplelatex}
\[
 - {\displaystyle \frac {k\,(-1)^{k}\,\mathrm{binomial}(n, \,k)}{
n}} 
\;.
\]
\end{maplelatex}

\end{maplegroup}
\noindent
We can increase the level of user information by the command

\begin{maplegroup}
\begin{mapleinput}
\mapleinline{active}{1d}{infolevel[sum]:=3:}{%
}
\end{mapleinput}

\end{maplegroup}
\noindent
Let's try yo prove the statement

\[
{\displaystyle \sum _{k=1}^{\infty }} \,{\displaystyle \frac {(-1
)^{k + 1}\,(4\,k + 1)\,(2\,k)\mathrm{!}}{k\mathrm{!}\,4^{k}\,(2
\,k - 1)\,(k + 1)\mathrm{!}}} =1
\]
that was posed in SIAM Review 36, 1994, Problem 94-2 {\rm \cite{SIAM}}.
We compute an antidifference

\begin{maplegroup}
\begin{mapleinput}

\mapleinline{active}{1d}{summand:=(-1)^(k+1)*(4*k+1)*(2*k)!/(k!*4^k*(2*k-1)*(k+1)!):}{%
}

\mapleinline{active}{1d}{sum(summand,k);
}{%
}
\end{mapleinput}

\mapleresult
\begin{maplettyout}
sum/indefnew:   indefinite summation
sum/extgosper:   applying Gosper algorithm to a(   k   ):=  
(-1)^(k+1)*(4*k+1)*(2*k)!/k!/(4^k)/(2*k-1)/(k+1)!
sum/gospernew:   a(   k   )/a(   k   -1):=   
-1/2*(4*k+1)/(4*k-3)/(k+1)*(2*k-3)
sum/gospernew:   Gosper's algorithm applicable
sum/gospernew:   p:=   4*k+1
sum/gospernew:   q:=   -2*k+3
sum/gospernew:   r:=   2*k+2
sum/gospernew:   degreebound:=   0
sum/gospernew:   solving equations to find f
sum/gospernew:   Gosper's algorithm successful
sum/gospernew:   f:=   -1
sum/indefnew:   indefinite summation finished
\end{maplettyout}

\begin{maplelatex}
\[
 - 2\,{\displaystyle \frac {(k + 1)\,(-1)^{(k + 1)}\,(2\,k)
\mathrm{!}}{k\mathrm{!}\,4^{k}\,(2\,k - 1)\,(k + 1)\mathrm{!}}} 
\;,
\]
\end{maplelatex}

\end{maplegroup}
\noindent
with success. Taking the limit as $n\rightarrow\infty$, one gets therefore

\begin{maplegroup}
\begin{mapleinput}

\mapleinline{active}{1d}{sum(summand,k=1..infinity);
}{%
}
\end{mapleinput}

\mapleresult
\begin{maplettyout}
sum/infinite:   infinite summation
\end{maplettyout}

\begin{maplelatex}
\[
1
\;.
\]
\end{maplelatex}

\end{maplegroup}
\noindent
Moreover, from the user information we see that {\sl Gosper's algorithm}
is applied.

If $a_k$ is a hypergeometric term, i.e., if%
\footnote{$\Q(k)$: rational functions over $\Q$.}
\[
\frac{a_{k+1}}{a_{k}}=\frac{b_k}{c_k}
\in\Q(k)
\;,
\]
then Gosper's algorithm {\sl decides} whether or not the antidifference 
$s_k$ is a hypergeometric term, and computes it in the affirmative case.

In detail: Given%
\footnote{$\Q[k]$: polynomials over $\Q$.}
\[
\frac{a_{k+1}}{a_{k}}=\frac{b_k}{c_k}
\;,\quad\quad
b_k, c_k \in\Q[k]
\;,
\]
a representation
\[
\frac{b_k}{c_k}
=\frac{p_{k+1}}{p_{k}}\,\frac{q_{k+1}}{r_{k+1}}
\;,\quad\quad
p_k, q_k, r_k \in\Q[k]
\]
is computed for which
\[
\gcd\:(q_k,r_{k+j})=1
\quad\quad\mbox{for all}\;j\in\N_0
\;.
\]
This can be done by a resultant computation \cite{Gosper} or by
rational factorization (\cite{Koornwinder}, \cite{MW}).

The essential fact is then: $f_k$, defined by
\[
s_{k}=\frac{r_{k}}{p_{k}}\,f_{k-1}\,a_{k}
\]
is rational, and the above gcd-condition yields even $f_k\in \Q[k]$.
$f_k$ satisfies the inhomogeneous recurrence equation
\[
p_k
=
q_{k+1}\,f_k-r_k\,f_{k-1}
\;.
\]
After calculating the degree of $f_k$, it is pure linear algebra to compute
$f_k$. The output of the procedure is either
\[
\displaystyle{s_{k}=\frac{r_{k}}{p_{k}}\,f_{k-1}\,a_{k}}
\]
or the statement {\sl ``There is no elementary (= hypergeometric term)
antidifference''}.

In the book \cite{Koepf}, many algorithms that are connected with Gosper's,
are treated in detail and Maple implementations are given.

After loading {\tt `hsum.mpl`},%
\footnote{The packages {\tt `hsum.mpl`} and {\tt `qsum.mpl`} can be
obtained from the URL {\tt www.imn.htwk-leipzig.de/\~{}koepf/research.html}.}

\begin{maplegroup}

\begin{mapleinput}
\mapleinline{active}{1d}{read(`hsum.mpl`);}{%
}
\end{mapleinput}

\mapleresult
\begin{maplelatex}
\[
\mathit{Copyright\ 1998\ \ Wolfram\ Koepf,\ Konrad-Zuse-Zentrum\ 
Berlin}
\]
\end{maplelatex}

\end{maplegroup}
\noindent
we can repeat the above calculation by the command

\begin{maplegroup}
\begin{mapleinput}

\mapleinline{active}{1d}{gosper((-1)^(k+1)*(4*k+1)*(2*k)!/(k!*4^k*(2*k-1)*(k+1)!),k);
}{%
}
\end{mapleinput}

\mapleresult
\begin{maplelatex}
\[
 - 2\,{\displaystyle \frac {(k + 1)\,(-1)^{(k + 1)}\,(2\,k)
\mathrm{!}}{k\mathrm{!}\,4^{k}\,(2\,k - 1)\,(k + 1)\mathrm{!}}} 
\;.
\]
\end{maplelatex}

\end{maplegroup}
\noindent
The computation

\begin{maplegroup}
\begin{mapleinput}
\mapleinline{active}{1d}{gosper(1/k,k);}{%
}
\end{mapleinput}

\mapleresult
\begin{maplettyout}
Error, (in gosper) no hypergeometric term antidifference exists
\end{maplettyout}

\end{maplegroup}
\begin{maplegroup}
\begin{mapleinput}
\end{mapleinput}

\end{maplegroup}
\noindent
is not worthless at all: It proves that the harmonic numbers
\[
H_n:=\sum_{k=1}^n\frac{1}{k}
\]
do {\sl not} constitute a hypergeometric term, corresponding to the
fact that the logarithmic function
\[
\ln x=\int_1^x\frac{1}{t}\,dt
\]
cannot be written in terms of exponentials.

Zeilberger's algorithms is an extension of Gosper's for definite sums.
It generates, e.g., the right-hand sides of the identities
\[
\sum_{k=0}^n {{n}\choose{k}}=
2^n
\;,
\]
\[
\sum_{k=0}^n {{n}\choose{k}}^2=
{\frac {\left (2\,n\right )!}{\left (n!\right )^{2}}}
\]
by the commands

\begin{maplegroup}

\begin{mapleinput}
\mapleinline{active}{1d}{closedform(binomial(n,k),k,n);}{%
}
\end{mapleinput}

\mapleresult
\begin{maplelatex}
\[
2^{n}
\]
\end{maplelatex}

\end{maplegroup}
\begin{maplegroup}
\begin{mapleinput}
\mapleinline{active}{1d}{closedform(binomial(n,k)^2,k,n);}{%
}
\end{mapleinput}

\mapleresult
\begin{maplelatex}
\[
{\displaystyle \frac {(2\,n)\mathrm{!}}{(n\mathrm{!})^{2}}} 
\;.
\]
\end{maplelatex}

\end{maplegroup}
\noindent
Here are the details:
If $F(n,k)$ is a hypergeometric term w.r.t.\ $n$ and $k$, i.e.
\[
\frac{F(n+1,k)}{F(n,k)}
\quad
\mbox{and}
\quad
\frac{F(n,k+1)}{F(n,k)}
\in\Q(n,k)
\;,
\]
then Zeilberger's algorithm generates a {\sl holonomic recurrence equation}
for
\[
s_n:=\sum_{k\in\Z} F(n,k)
\;.
\]
This is performed by starting with $J=1$ and iterating if necessary: Set
\[
a_k:=F(n,k)+\sum_{j=1}^J \sigma_j(n)\,F(n+j,k)
\]
with as yet undetermined variables $\sigma_j$.
Apply Gosper's algorithm to $a_k$.
In the last step, solve the linear system at the same time for
the coefficients of $f_k$ and the variables $\sigma_j\;(j=1,\ldots,J)$.
In the affirmative case, this yields
\[
G(n,k+1)-G(n,k)=a_k
\;.
\]
Output:
By summation one gets:
\[
s_n+\sum_{j=1}^J \sigma_j(n)\,s_{n+j}
=0
\;.
\]
We would like to point out that the most time consuming part of Zeilberger's
algorithm is its last step which is to solve a linear system. This linear
system, however, often has many variables, and its coefficients are
polynomials or rational functions. Here, an efficient implementation of
linear algebra is important. Furthermore, the resulting recurrence equation
usually needs factored coefficients because otherwise the results look
unnecessarily complicated. We will see such a situation soon.

We give some examples:
Each of the following series represents the {\sl Legendre polynomials}:
\begin{eqnarray*}
P_n(x)&=&
\sum_{k=0}^n {n\choose k}\,{-n-1\choose k}\left(\frac{1-x}{2}\right)^k
\\&=&
\hypergeom{2}{1}{-n,n+1}{1}{\frac{1-x}{2}}
\\&=&
\frac{1}{2^n}\,\sum_{k=0}^n
{n\choose k}^2\,(x-1)^{n-k}\,(x+1)^k
\\&=&
\left( \frac{1-x}{2} \right)^n\;\hypergeom{2}{1}{-n,-n}{1}{\frac{1+x}{1-x}}
\\&=&
\frac{1}{2^n}\,\sum_{k=0}^{\lfloor n/2\rfloor}
(-1)^k\,{n\choose k}\,{2n-2k\choose n}\
\,x^{n-2k}
\\&=&
{2n\choose n}\,\left(\frac{x}{2}\right)^n
\hypergeom{2}{1}{-n/2,-n/2+1/2}{-n+1/2}{\frac{1}{x^2}}
\\[8mm]&=&
x^n\;\hypergeom{2}{1}{-n/2,-n/2+1/2}{1}{1-\frac{1}{x^2}}
\;.
\end{eqnarray*}
Again, you see, that functions come in quite different disguises.
How can we show that these systems define the same family of functions?
Zeilberger's paradigm is to show that they satisfy the same recurrence
equation, then it is sufficient to check a finite number of initial values.

Here are the recurrence equations for the different sums:

\begin{maplegroup}

\begin{mapleinput}
\mapleinline{active}{1d}{P:='P':}{%
}
\end{mapleinput}

\end{maplegroup}
\begin{maplegroup}
\begin{mapleinput}

\mapleinline{active}{1d}{sumrecursion(binomial(n,k)*binomial(-n-1,k)*((1-x)/2)^k,k,P(n));
}{%
}
\end{mapleinput}

\mapleresult
\begin{maplelatex}
\[
(n + 2)\,\mathrm{P}(n + 2) - (2\,n + 3)\,x\,\mathrm{P}(n + 1) + (
n + 1)\,\mathrm{P}(n)=0
\]
\end{maplelatex}

\end{maplegroup}
\begin{maplegroup}
\begin{mapleinput}

\mapleinline{active}{1d}{sumrecursion(1/2^n*binomial(n,k)^2*(x-1)^(n-k)*(x+1)^k,k,P(n));
}{%
}
\end{mapleinput}

\mapleresult
\begin{maplelatex}
\[
(n + 2)\,\mathrm{P}(n + 2) - (2\,n + 3)\,x\,\mathrm{P}(n + 1) + (
n + 1)\,\mathrm{P}(n)=0
\]
\end{maplelatex}

\end{maplegroup}
\begin{maplegroup}
\begin{mapleinput}
\mapleinline{active}{1d}{sumrecursion(1/2^n*(-1)^k*binomial(n,k)*}{%
}
\mapleinline{active}{1d}{binomial(2*n-2*k,n)*x^(n-2*k),k,P(n));}{%
}
\end{mapleinput}

\mapleresult
\begin{maplelatex}
\[
(n + 2)\,\mathrm{P}(n + 2) - (2\,n + 3)\,x\,\mathrm{P}(n + 1) + (
n + 1)\,\mathrm{P}(n)=0
\]
\end{maplelatex}

\end{maplegroup}
\begin{maplegroup}
\begin{mapleinput}

\mapleinline{active}{1d}{sumrecursion(x^n*hyperterm([-n/2,-n/2+1/2],[1],1-1/x^2,k),k,P(n));
}{%
}
\end{mapleinput}

\mapleresult
\begin{maplelatex}
\[
(n + 2)\,\mathrm{P}(n + 2) - (2\,n + 3)\,x\,\mathrm{P}(n + 1) + (
n + 1)\,\mathrm{P}(n)=0
\]
\end{maplelatex}

\end{maplegroup}
\noindent
We omit the computation of the initial values. 

The {\tt Sumtohyper}
procedure of the {\tt hsum} package is slightly more efficient than
{\tt `convert/hypergeom`} by
converting a series into hypergeometric notation.

\begin{maplegroup}

\begin{mapleinput}

\mapleinline{active}{1d}{Sumtohyper(binomial(n,k)*binomial(-n-1,k)*((1-x)/2)^k,k);
}{%
}
\end{mapleinput}

\mapleresult
\begin{maplelatex}
\[
\mathrm{Hypergeom}([n + 1, \, - n], \,[1], \,{\displaystyle 
\frac {1}{2}}  - {\displaystyle \frac {1}{2}} \,x)
\]
\end{maplelatex}

\end{maplegroup}
\begin{maplegroup}
\begin{mapleinput}

\mapleinline{active}{1d}{Sumtohyper(1/2^n*binomial(n,k)^2*(x-1)^(n-k)*(x+1)^k,k);
}{%
}
\end{mapleinput}

\mapleresult
\begin{maplelatex}
\[
({\displaystyle \frac {1}{2}} \,x - {\displaystyle \frac {1}{2}} 
)^{n}\,\mathrm{Hypergeom}([ - n, \, - n], \,[1], \,
{\displaystyle \frac {x + 1}{x - 1}} )
\]
\end{maplelatex}

\end{maplegroup}
\begin{maplegroup}
\begin{mapleinput}

\mapleinline{active}{1d}{Sumtohyper(1/2^n*(-1)^k*binomial(n,k)*}{%
}
\mapleinline{active}{1d}{binomial(2*n-2*k,n)*x^(n-2*k),k);
}{%
}
\end{mapleinput}

\mapleresult
\begin{maplelatex}
\[
2^{( - n)}\,\mathrm{binomial}(2\,n, \,n)\,x^{n}\,\mathrm{
Hypergeom}([ - {\displaystyle \frac {1}{2}} \,n + {\displaystyle 
\frac {1}{2}} , \, - {\displaystyle \frac {1}{2}} \,n], \,[ - n
 + {\displaystyle \frac {1}{2}} ], \,{\displaystyle \frac {1}{x^{
2}}} )
\;.
\]
\end{maplelatex}

\end{maplegroup}
\noindent
The above computations show that all the given representations of the
Legendre polynomials agree.

To give a more advanced example of an application of {\tt Sumtohyper},
we compute the hypergeometric representation of the difference
$P_{n+1}(x)-P_n(x)$ of successive Legendre polynomials:

\begin{maplegroup}
\begin{mapleinput}

\mapleinline{active}{1d}{legendreterm:=binomial(n,k)*binomial(-n-1,k)*((1-x)/2)^k;
}{%
}
\end{mapleinput}

\mapleresult
\begin{maplelatex}
\[
\mathit{legendreterm} := \mathrm{binomial}(n, \,k)\,\mathrm{
binomial}( - n - 1, \,k)\,({\displaystyle \frac {1}{2}}  - 
{\displaystyle \frac {1}{2}} \,x)^{k}
\]
\end{maplelatex}

\end{maplegroup}

\begin{maplegroup}
\begin{mapleinput}

\mapleinline{active}{1d}{Sumtohyper(subs(n=n+1,legendreterm)-legendreterm,k);
}{%
}
\end{mapleinput}

\mapleresult
\begin{maplelatex}
\[
(x + x\,n - 1 - n)\,\mathrm{Hypergeom}([ - n, \,n + 2], \,[2], \,
{\displaystyle \frac {1}{2}}  - {\displaystyle \frac {1}{2}} \,x)
\;.
\]
\end{maplelatex}

\end{maplegroup}
\noindent
We give more examples of how Zeilberger's algorithm can be applied in
rather different situations.

The following recurrence equation of the {\sl Ap\'ery numbers}
\[
A_n:=\sum_{k=0}^n {n \choose k}^2\,{n+k \choose k}^2
\]
was an essential tool in Ap\'ery's proof of the irrationality of
\[
\zeta (3)=\sum_{j=1}^\infty \frac{1}{j^3}
\;:
\]

\begin{maplegroup}
\begin{mapleinput}

\mapleinline{active}{1d}{sumrecursion(binomial(n,k)^2*binomial(n+k,k)^2,k,A(n));
}{%
}
\end{mapleinput}

\mapleresult
\begin{maplelatex}
\[
{(n + 2)^{3}\,\mathrm{A}(n + 2) - (3 + 2\,n)\,(17\,n^{2}
 + 51\,n + 39)\,\mathrm{A}(n + 1) + (n + 1)^{3}\,\mathrm{A}(n)}
  =0
\;.
\]
\end{maplelatex}

\end{maplegroup}
\begin{maplegroup}
\begin{mapleinput}
\end{mapleinput}

\end{maplegroup}
\noindent
{\sl Dougall's identity}
\[
_7 F_6\left.
\!\!
\left(
\!\!\!\!
\begin{array}{c}
\multicolumn{1}{c}{\begin{array}{c}
a,1+\frac{a}{2},b,c,d,1+2a-b-c-d+n,-n
\end{array}}\\[5mm]
\multicolumn{1}{c}{\begin{array}{c}
\frac{a}{2},1\!+\!a\!-\!b,1\!+\!a\!-\!c,1\!+\!a\!-\!d,b\!+\!c\!+\!d\!-\!a\!-\!n,
1\!+\!a\!+\!n
            \end{array}}\end{array}
\!\!\!\!
\right| 1\right)
=
\]
\begin{equation}
\frac{(1+a)_n\,(1+a-b-c)_{n}\,(1+a-b-d)_n\,(1+a-c-d)_n}
{(1+a-b)_n(1+a-c)_n\,(1+a-d)_{n}\,(1+a-b-c-d)_n}
\label{eq:Dougall}
\end{equation}
is proven by

\begin{maplegroup}

\begin{mapleinput}

\mapleinline{active}{1d}{sumrecursion(hyperterm([a,1+a/2,b,c,d,1+2*a-b-c-d+n,-n],}{%
}
\mapleinline{active}{1d}{[a/2,1+a-b,1+a-c,1+a-d,b+c+d-a-n,1+a+n],1,k),k,S(n));
}{%
}
\end{mapleinput}

\mapleresult
\begin{maplelatex}
\begin{eqnarray*}
\lefteqn{ - (a - d + n + 1)\,(n + 1 + a - c)\,(n + 1 - b + a)\,(
 - b - c - d + a + n + 1)} \\
 & & \mathrm{S}(n + 1)\mbox{} + (1 + a + n)\,(n + 1 - c + a - d)
\,(n + 1 - b + a - d) \\
 & & (n - c + a + 1 - b)\,\mathrm{S}(n)=0
\;.
\mbox{\hspace{170pt}}
\end{eqnarray*}
\end{maplelatex}

\end{maplegroup}
\noindent
From this result, the right-hand side (\ref{eq:Dougall}) of Dougall's
identity can be read off directly. The complete computation is performed by

\begin{maplegroup}
\begin{mapleinput}

\mapleinline{active}{1d}{closedform(hyperterm([a,1+a/2,b,c,d,1+2*a-b-c-d+n,-n],}{%
}
\mapleinline{active}{1d}{[a/2,1+a-b,1+a-c,1+a-d,b+c+d-a-n,1+a+n],1,k),k,n);
}{%
}
\end{mapleinput}

\mapleresult
\begin{maplelatex}
\begin{eqnarray*}
\lefteqn{\mathrm{pochhammer}(a + 1, \,n)\,\mathrm{pochhammer}( - 
d - c + 1 + a, \,n)} \\
 & & \mathrm{pochhammer}( - d + 1 - b + a, \,n)\,\mathrm{
pochhammer}(a + 1 - b - c, \,n) \left/ {\vrule 
height0.37em width0em depth0.37em} \right. \!  \! ( \\
 & & \mathrm{pochhammer}(1 + a - d, \,n)\,\mathrm{pochhammer}(1
 + a - c, \,n) \\
 & & \mathrm{pochhammer}(1 + a - b, \,n)\,\mathrm{pochhammer}(a
 - d + 1 - b - c, \,n))
\;.
\end{eqnarray*}
\end{maplelatex}

\end{maplegroup}
\noindent
Notice how important rational factorization is for such examples!

The {\sl Wilson polynomials} have the representation%
\footnote{Sometimes a different standardization is used. But this is not 
essential.}
\[
W_n(x)
=
\hypergeom{4}{3}{-n,a\!+\!b\!+\!c\!+\!d\!+\!n\!-\!1,a\!-\!x,a\!+\!x}
{a+b,a+c,a+d}{1}
\;.
\]
They include all classical systems like the Jacobi and Hahn polynomials.
We get

\begin{maplegroup}

\begin{mapleinput}

\mapleinline{active}{1d}{sumrecursion(hyperterm([-n,a+b+c+d+n-1,a-x,a+x],}{%
}
\mapleinline{active}{1d}{[a+b,a+c,a+d],1,k),k,W(n));
}{%
}
\end{mapleinput}

\mapleresult
\begin{maplelatex}
\begin{eqnarray*}
\lefteqn{(d + a + n + 1)\,(n + 1 + a + c)\,(n + b + 1 + a)\,(a + 
2\,n + c + b + d)} \\
 & & (a + b + c + d + n)\,\mathrm{W}(n + 2)\mbox{} - (2\,n + 1 + 
a + b + c + d)(8\,c\,d\,n \\
 & & \mbox{} + 8\,b\,d\,n + b\,d\,a^{2} + 8\,b\,c\,n + 3\,b^{2}\,
n + b\,c\,a^{2} + 2\,b^{2}\,d + 6\,a\,n\,b\,c \\
 & & \mbox{} + 7\,d\,n^{2} + 6\,b\,n\,c\,d + 3\,d^{2}\,n + 2\,c\,
d^{2} + 2\,b\,d^{2} + 7\,c\,n^{2} + 3\,c^{2}\,n \\
 & & \mbox{} + 2\,c^{2}\,d + 2\,b\,c^{2} + 7\,b\,n^{2} + 2\,b^{2}
\,c + 7\,a\,n^{2} + 8\,a\,d\,n + 4\,n^{3}\,a \\
 & & \mbox{} + 2\,d^{2}\,n^{2} + 4\,d\,n^{3} + 2\,c^{2}\,n^{2} + 
4\,c\,n^{3} + 2\,b^{2}\,n^{2} + 4\,b\,n^{3} + 6\,a\,n\,b\,d \\
 & & \mbox{} + b^{2} + c^{2} + 8\,a\,c\,n + d^{2} + 8\,a\,b\,n + 
4\,a\,b\,c + 2\,a\,b + 6\,a\,n\,c\,d \\
 & & \mbox{} + b\,c\,d^{2} + b\,c^{2}\,d + a\,b^{2}\,d + 2\,a\,b
^{2}\,n + a\,c\,d^{2} + 2\,a\,c^{2}\,n + a\,b\,d^{2} \\
 & & \mbox{} + a\,b^{2}\,c + a\,b\,c^{2} + b^{2}\,c\,d + 4\,a\,b
\,c\,d + a\,c^{2}\,d + 3\,a^{2}\,n + 6\,b\,n^{2}\,a \\
 & & \mbox{} + 6\,b\,n^{2}\,d + 6\,b\,n^{2}\,c + 2\,b^{2}\,n\,d
 + 2\,b^{2}\,n\,c + 4\,n^{3} + 2\,n^{4} + 2\,c^{2}\,n\,b \\
 & & \mbox{} + 2\,a\,d^{2} + 2\,a\,c^{2} + 2\,a\,b^{2} + 6\,c\,n
^{2}\,a + 6\,c\,n^{2}\,d + 2\,c^{2}\,n\,d + 2\,a^{2}\,c\,n \\
 & & \mbox{} + 2\,a^{2}\,b\,n + 2\,x^{2}\,b\,d + 2\,x^{2}\,c\,d
 + 2\,x^{2}\,a\,d + x^{2}\,a^{2} + 2\,x^{2}\,b\,c \\
 & & \mbox{} + 2\,x^{2}\,a\,c + x^{2}\,d^{2} + 2\,x^{2}\,a\,b + x
^{2}\,b^{2} + x^{2}\,c^{2} + 4\,b\,c\,d + 2\,d^{2}\,n\,c \\
 & & \mbox{} + 2\,d^{2}\,n\,b + 2\,a^{2}\,d + 6\,d\,n^{2}\,a + 4
\,a\,b\,d + 2\,n^{2} + 2\,a^{2}\,b + 2\,a^{2}\,c \\
 & & \mbox{} + 2\,b\,c + 4\,a\,c\,d + 2\,a\,d^{2}\,n + a^{2} + 2
\,a\,c + 3\,b\,n + 2\,b\,d + 2\,c\,d \\
 & & \mbox{} + 3\,d\,n + 2\,a\,d + 2\,a^{2}\,d\,n + 3\,a\,n + 4\,
a\,x^{2}\,n + 4\,b\,x^{2}\,n + 4\,c\,x^{2}\,n \\
 & & \mbox{} + 4\,d\,x^{2}\,n + 4\,x^{2}\,n^{2} + 2\,a\,x^{2} + 2
\,b\,x^{2} + 2\,c\,x^{2} + 2\,d\,x^{2} + 4\,x^{2}\,n \\
 & & \mbox{} + 3\,c\,n + 2\,a^{2}\,n^{2} + c\,d\,a^{2})\mathrm{W}
(n + 1)\mbox{} + (n + 1)\,(n + d + c)\,(n + b + d) \\
 & & (n + b + c)\,(d + 2\,n + a + b + c + 2)\,\mathrm{W}(n)=0
\;,
\end{eqnarray*}
\end{maplelatex}

\end{maplegroup}
\noindent
a recurrence equation for $W_n(x)$ which, however, is rather complicated
since the middle coefficient admits no rational factorization.

One knows from the theory that the recurrence equation has a
special form which can be found by the command {\tt Sumrecursion}:

\begin{maplegroup}
\begin{mapleinput}

\mapleinline{active}{1d}{Sumrecursion(hyperterm([-n,n+a+b+c+d-1,a+x,a-x],}{%
}
\mapleinline{active}{1d}{[a+b,a+c,a+d],1,k),k,W(n,x));
}{%
}
\end{mapleinput}

\mapleresult
\begin{maplelatex}
\begin{eqnarray*}
\lefteqn{(x - a)\,(a + x)\,\mathrm{W}(n, \,x)=((a + d + n)\,(n +
a + c)\,(n + b + a)} \\
 & & (a + b + c + d + n - 1)\,\mathrm{W}(n + 1, \,x)) \left/
{\vrule height0.37em width0em depth0.37em} \right. \!  \! ( \\
 & & (a + 2\,n + c + b + d)\,(2\,n + b + a + c + d - 1))\mbox{}
 - ( \\
 & & {\displaystyle \frac {(a + d + n)\,(n + a + c)\,(n + b + a)
\,(a + b + c + d + n - 1)}{(a + 2\,n + c + b + d)\,(2\,n + b + a
 + c + d - 1)}}  \\
 & & \mbox{} + {\displaystyle \frac {(n + d + c - 1)\,(n + b + d
 - 1)\,(n + b - 1 + c)\,n}{(2\,n + b + a + c + d - 1)\,(2\,n - 2
 + b + d + a + c)}} ) \\
 & & \mathrm{W}(n, \,x)\mbox{} +  \\
 & & {\displaystyle \frac {(n + d + c - 1)\,(n + b + d - 1)\,(n
 + b - 1 + c)\,n\,\mathrm{W}(n - 1, \,x)}{(2\,n + b + a + c + d
 - 1)\,(2\,n - 2 + b + d + a + c)}}
\;.
\end{eqnarray*}
\end{maplelatex}

\end{maplegroup}
\noindent
A similar recurrence equation exists w.r.t.\ $x$:

\begin{maplegroup}
\begin{mapleinput}

\mapleinline{active}{1d}{Sumrecursion(hyperterm([-n,n+a+b+c+d-1,a+x,a-x],}{%
}
\mapleinline{active}{1d}{[a+b,a+c,a+d],1,k),k,W(x,n));
}{%
}
\end{mapleinput}

\mapleresult
\begin{maplelatex}
\begin{eqnarray*}
\lefteqn{(a + b + c + d + n - 1)\,n\,\mathrm{W}(x, \,n)=} \\
 & & {\displaystyle \frac {1}{2}} \,{\displaystyle \frac {(x + d)
\,(x + c)\,(x + b)\,(a + x)\,\mathrm{W}(x + 1, \,n)}{(2\,x + 1)\,
x}}  -  \\
 & & ({\displaystyle \frac {1}{2}} \,{\displaystyle \frac {(x + d
)\,(x + c)\,(x + b)\,(a + x)}{(2\,x + 1)\,x}}  + {\displaystyle 
\frac {1}{2}} \,{\displaystyle \frac {(x - d)\,(x - c)\,(x - b)\,
(x - a)}{( - 1 + 2\,x)\,x}} ) \\
 & & \mathrm{W}(x, \,n)\mbox{} + {\displaystyle \frac {1}{2}} \,
{\displaystyle \frac {(x - d)\,(x - c)\,(x - b)\,(x - a)\,
\mathrm{W}(x - 1, \,n)}{( - 1 + 2\,x)\,x}} 
\;.
\end{eqnarray*}
\end{maplelatex}

\end{maplegroup}
\noindent
{\sl Clausen's formula}
\[
\hypergeom{2}{1}{a,b}{a\!+\!b\!+\!1/2}{x}^2
=
\hypergeom{3}{2}{2a,2b,a+b}{a\!+\!b\!+\!1/2,2a\!+\!2b}{x}
\]
gives the cases when the square of a $_2 F_1$ function is a $_3 F_2$. 
The right-hand side is deduced from the left-hand side by

\begin{maplegroup}

\begin{mapleinput}

\mapleinline{active}{1d}{sumrecursion(hyperterm([a,b],[a+b+1/2],x,j)*}{%
}
\mapleinline{active}{1d}{hyperterm([a,b],[a+b+1/2],x,k-j),j,C(k));
}{%
}
\end{mapleinput}

\mapleresult
\begin{maplelatex}
\begin{eqnarray*}
\lefteqn{ - (k + 1)\,(2\,a + 1 + 2\,b + 2\,k)\,(2\,a + 2\,b + k)
\,\mathrm{C}(k + 1)} \\
 & & \mbox{} + 2\,x\,(k + 2\,b)\,(k + 2\,a)\,(a + b + k)\,
\mathrm{C}(k)=0
\;,
\mbox{\hspace{4pt}}
\end{eqnarray*}
\end{maplelatex}

\end{maplegroup}
\noindent
computing the coefficient of the Cauchy product. The resulting hypergeometric
term can be obtained in one step by%
\footnote{The {\tt Closedform} procedure differs from the 
{\tt closedform} procedure in that the hypergeomtric term is not
evaluated.}

\begin{maplegroup}
\begin{mapleinput}

\mapleinline{active}{1d}{Closedform(hyperterm([a,b],[a+b+1/2],x,j)*}{%
}
\mapleinline{active}{1d}{hyperterm([a,b],[a+b+1/2],x,k-j),j,k);
}{%
}
\end{mapleinput}

\mapleresult
\begin{maplelatex}
\[
\mathrm{Hyperterm}([2\,b, \,2\,a, \,b + a], \,[a + b + 
{\displaystyle \frac {1}{2}} , \,2\,a + 2\,b], \,x, \,k)
\]
\end{maplelatex}

\end{maplegroup}
\begin{maplegroup}
\begin{mapleinput}
\end{mapleinput}

\end{maplegroup}

\begin{maplegroup}
\noindent
The computation of a specific {\sl Feynman diagram} \cite{FT}
yields the representation
\[
V(\alpha,\beta,\gamma)
=
(-1)^{\alpha+\beta+\gamma}\;\cdot
\frac{\Gamma(\alpha\!+\!\beta\!+\!\gamma\!-\!d/2)
\Gamma(d/2\!-\!\gamma)\Gamma(\alpha\!+\!\gamma\!-\!d/2)
\Gamma(\beta\!+\!\gamma\!-\!d/2)}
{\Gamma(\alpha)\Gamma(\beta)\Gamma(d/2)\Gamma(\alpha+\beta+2\gamma-d)
(m^2)^{\alpha+\beta+\gamma-d}}
\]
\[
\cdot
\; _2 F_1\left.
\!\!
\left(
\!\!\!\!
\begin{array}{c}
\multicolumn{1}{c}{\begin{array}{cc} \alpha+\beta+\gamma-d
\;, & \alpha+\gamma-d/2 \end{array}}\\[1mm]
\multicolumn{1}{c}{ \alpha+\beta+2\gamma-d}
            \end{array}
\!\!\!\!
\right| z\right)
\;.
\]
Since one is interested to compute this function for $\alpha,\beta,\gamma\in\N$,
and since the computation is easy for $\alpha,\beta,\gamma\in\{0,1\}$,
recurrence equations w.r.t.\ these variables can be used. Here is one w.r.t.\ 
$\beta$:

\begin{mapleinput}
\end{mapleinput}

\end{maplegroup}
\begin{maplegroup}

\begin{mapleinput}
\mapleinline{active}{1d}{sumrecursion((-1)^(alpha+beta+gamma)*}{%
}
\mapleinline{active}{1d}{GAMMA(alpha+beta+gamma-d/2)*GAMMA(d/2-gamma)*}{%
}
\mapleinline{active}{1d}{
GAMMA(alpha+gamma-d/2)*GAMMA(beta+gamma-d/2)/}{%
}
\mapleinline{active}{1d}{
(GAMMA(alpha)*GAMMA(beta)*GAMMA(d/2)*}{%
}
\mapleinline{active}{1d}{
GAMMA(alpha+beta+2*gamma-d)*(m^2)^}{%
}
\mapleinline{active}{1d}{
(alpha+beta+gamma-d))*}{%
}
\mapleinline{active}{1d}{
hyperterm([alpha+beta+gamma-d,}{%
}
\mapleinline{active}{1d}{
alpha+gamma-d/2],[alpha+beta+2*gamma-d],z,k),k,V(beta));}{%
}
\end{mapleinput}

\mapleresult
\begin{maplelatex}
\begin{eqnarray*}
\lefteqn{8\,\beta \,m^{4}\,(\beta  + 1)\,(\alpha  + \beta  + 1 + 
\gamma  - d)\,z\,\mathrm{V}(\beta  + 2) + 2\beta \,m^{2}} \\
 & & (4\,\gamma  + 2\,z\,\beta  + 2\,z - z\,d + 2\,\alpha  + 2\,
\beta  - 2\,d)\,(2\,\alpha  + 2\,\beta  + 2 + 2\,\gamma  - d) \\
 & & \mathrm{V}(\beta  + 1)\mbox{} +  \\
 & & (2\,\alpha  + 2\,\beta  + 2\,\gamma  - d)\,(2\,\gamma  - d
 + 2\,\beta )\,(2\,\alpha  + 2\,\beta  + 2 + 2\,\gamma  - d)\,
\mathrm{V}(\beta ) \\
 & & =0
\;.
\end{eqnarray*}
\end{maplelatex}

\end{maplegroup}
\noindent
Similarly, one obtains recurrence equations w.r.t.\ $\alpha$ and $\gamma$.

In some instances, Zeilberger's algorithm does not find the recurrence
equation of lowest order. Assume, e.g., we want to deduce the
right-hand side from the left-hand side of the identity
\[
{\displaystyle \sum _{k=0}^{n}} \,(-1)^{k}\,{n \choose k}\,{3k \choose n}=
(-3)^{n}
\;.
\]
Therefore, by Zeilberger's algorithm we compute a recurrence equation 

\begin{maplegroup}
\begin{mapleinput}

\mapleinline{active}{1d}{RE:=sumrecursion((-1)^k*binomial(n,k)*binomial(3*k,n),k,S(n));
}{%
}
\end{mapleinput}

\mapleresult
\begin{maplelatex}
\[
\mathit{RE} := 2\,(3 + 2\,n)\,\mathrm{S}(n + 2) + 3\,(7 + 5\,n)\,
\mathrm{S}(n + 1) + 9\,(n + 1)\,\mathrm{S}(n)=0
\]
\end{maplelatex}

\end{maplegroup}
\noindent
and apply {\sl Petkov\v{s}ek's algorithm}
to find its hypergeometric term solutions

\begin{maplegroup}
\begin{mapleinput}
\mapleinline{active}{1d}{rechyper(RE,S(n));}{%
}
\end{mapleinput}

\mapleresult
\begin{maplelatex}
\[
\{-3\}
\]
\end{maplelatex}

\end{maplegroup}
\noindent
which gives the term ratio $S_{n+1}/S_n$ of the resulting hypergeometric term
$S_n=(-3)^n$.

Here is another application of Petkov\v{s}ek's algorithm:

\begin{maplegroup}

\begin{mapleinput}
\mapleinline{active}{1d}{RE:=(n+4)*s(n+2)+s(n+1)-(n+1)*s(n)=0;}{%
}
\end{mapleinput}

\mapleresult
\begin{maplelatex}
\[
\mathit{RE} := (n + 4)\,\mathrm{s}(n + 2) + \mathrm{s}(n + 1) - (
n + 1)\,\mathrm{s}(n)=0
\]
\end{maplelatex}

\end{maplegroup}
\begin{maplegroup}
\begin{mapleinput}
\mapleinline{active}{1d}{rechyper(RE,s(n));}{%
}
\end{mapleinput}

\mapleresult
\begin{maplelatex}
\[
\{{\displaystyle \frac {n + 1}{n + 3}} , \, - {\displaystyle 
\frac {(5 + 2\,n)\,(n + 1)}{(3 + 2\,n)\,(n + 3)}} \}
\;.
\]
\end{maplelatex}

\end{maplegroup}
\noindent
If continuous variables are involved, one can also compute {\sl holonomic
differential equations} for sums by a Zeilberger type algorithm
which is implemented in the {\tt sumdiffeq} procedure.

We take some of the series representations of the Legendre polynomials to deduce
the corresponding differential equation:

\begin{maplegroup}

\begin{mapleinput}

\mapleinline{active}{1d}{sumdiffeq(binomial(n,k)*binomial(-n-1,k)*((1-x)/2)^k,k,P(x));
}{%
}
\end{mapleinput}

\mapleresult
\begin{maplelatex}
\[
 - ( - 1 + x)\,(x + 1)\,({\frac {\partial ^{2}}{\partial x^{2}}}
\,\mathrm{P}(x)) - 2\,x\,({\frac {\partial }{\partial x}}\,
\mathrm{P}(x)) + \mathrm{P}(x)\,n\,(n + 1)=0
\]
\end{maplelatex}

\end{maplegroup}
\begin{maplegroup}
\begin{mapleinput}

\mapleinline{active}{1d}{sumdiffeq(1/2^n*binomial(n,k)^2*(x-1)^(n-k)*(x+1)^k,k,P(x));
}{%
}
\end{mapleinput}

\mapleresult
\begin{maplelatex}
\[
 - ( - 1 + x)\,(x + 1)\,({\frac {\partial ^{2}}{\partial x^{2}}}
\,\mathrm{P}(x)) - 2\,x\,({\frac {\partial }{\partial x}}\,
\mathrm{P}(x)) + \mathrm{P}(x)\,n\,(n + 1)=0
\;.
\]
\end{maplelatex}

\end{maplegroup}
\noindent
To show a quadratic transformation like
\[
\hypergeom{2}{1}{a,b}{2b}{\frac{4x}{(1+x)^2}}
=
(1+x)^{2a}\cdot
\hypergeom{2}{1}{a,a-b+1/2}{b+1/2}{x^2}
\;,
\]
we prove that both sides satisfy the same differential equation,
and show that enough initial values agree.

\begin{maplegroup}

\begin{mapleinput}

\mapleinline{active}{1d}{sumdiffeq(hyperterm([a,b],[2*b],4*x/(1+x)^2,k),k,Q(x));
}{%
}
\end{mapleinput}

\mapleresult
\begin{maplelatex}
\begin{eqnarray*}
\lefteqn{ - x\,(x - 1)\,(1 + x)^{2}\,({\frac {\partial ^{2}}{
\partial x^{2}}}\,\mathrm{Q}(x)) + 2\,(1 + x)\,( - x^{2} + b\,x^{
2} - 2\,x\,a + b)\,({\frac {\partial }{\partial x}}\,\mathrm{Q}(x
))} \\
 & & \mbox{} + 4\,\mathrm{Q}(x)\,(x - 1)\,a\,b=0
\mbox{\hspace{198pt}}
\end{eqnarray*}
\end{maplelatex}

\end{maplegroup}
\begin{maplegroup}
\begin{mapleinput}

\mapleinline{active}{1d}{sumdiffeq((1+x)^(2*a)*hyperterm([a,a-b+1/2],[b+1/2],x^2,k),k,Q(x));
}{%
}
\end{mapleinput}

\mapleresult
\begin{maplelatex}
\begin{eqnarray*}
\lefteqn{ - x\,(x - 1)\,(1 + x)^{2}\,({\frac {\partial ^{2}}{
\partial x^{2}}}\,\mathrm{Q}(x)) + 2\,(1 + x)\,( - x^{2} + b\,x^{
2} - 2\,x\,a + b)\,({\frac {\partial }{\partial x}}\,\mathrm{Q}(x
))} \\
 & & \mbox{} + 4\,\mathrm{Q}(x)\,(x - 1)\,a\,b=0
\mbox{\hspace{198pt}}
\end{eqnarray*}
\end{maplelatex}

\end{maplegroup}
\begin{maplegroup}
\begin{mapleinput}

\mapleinline{active}{1d}{eval([hypergeom([a,b],[2*b],4*x/(1+x)^2)=}{%
}
\mapleinline{active}{1d}{(1+x)^(2*a)*hypergeom([a,a-b+1/2],[b+1/2],x^2),}{%
}
\mapleinline{active}{1d}{diff(hypergeom([a,b],[2*b],4*x/(1+x)^2)=}{%
}
\mapleinline{active}{1d}{(1+x)^(2*a)*hypergeom([a,a-b+1/2],[b+1/2],x^2),x)],x=0);
}{%
}
\end{mapleinput}

\mapleresult
\begin{maplelatex}
\[
[1=1, \,2\,a=2\,a]
\;.
\]
\end{maplelatex}

\end{maplegroup}
\noindent
On p.\ 258 in Ramanujan's second notebook one finds the identity
\[
\hypergeom{2}{1}{\frac{1}{3},\frac{2}{3}}{1}{1-\left(\frac{1-x}{1+2x}\right)^3}
=
(1+2x)\;\hypergeom{2}{1}{\frac{1}{3},\frac{2}{3}}{1}{x^3}
\;.
\]
With Garvan we can ask the question:
For which $A,B,C,a,b,c,d$ is
\[
\hypergeom{2}{1}{A,B}{C}{1-\left(\frac{1-x}{1+2x}\right)^3}
=
(1+2x)^d\;\hypergeom{2}{1}{a,b}{c}{x^3}
\;?
\]
A computation with Maple gives

\begin{maplegroup}

\begin{mapleinput}

\mapleinline{active}{1d}{first:=hyperterm([A,B],[C],1-((1-x)/(1+2*x))^3,k);
}{%
}
\end{mapleinput}

\mapleresult
\begin{maplelatex}
\[
\mathit{first} := {\displaystyle \frac {\mathrm{pochhammer}(A, \,
k)\,\mathrm{pochhammer}(B, \,k)\,(1 - {\displaystyle \frac {(1 - 
x)^{3}}{(1 + 2\,x)^{3}}} )^{k}}{\mathrm{pochhammer}(C, \,k)\,k
\mathrm{!}}} 
\]
\end{maplelatex}

\end{maplegroup}
\begin{maplegroup}
\begin{mapleinput}
\mapleinline{active}{1d}{second:=(2*x+1)^d*hyperterm([a,b],[c],x^3,k);
}{%
}
\end{mapleinput}

\mapleresult
\begin{maplelatex}
\[
\mathit{second} := {\displaystyle \frac {(1 + 2\,x)^{d}\,\mathrm{
pochhammer}(a, \,k)\,\mathrm{pochhammer}(b, \,k)\,(x^{3})^{k}}{
\mathrm{pochhammer}(c, \,k)\,k\mathrm{!}}} 
\]
\end{maplelatex}

\end{maplegroup}
\begin{maplegroup}
\begin{mapleinput}
\mapleinline{active}{1d}{DE1:=sumdiffeq(first,k,S(x));}{%
}
\end{mapleinput}

\mapleresult
\begin{maplelatex}
\begin{eqnarray*}
\lefteqn{\mathit{DE1} := x\,(x - 1)\,(1 + x + x^{2})\,(1 + 2\,x)
^{2}\,({\frac {\partial ^{2}}{\partial x^{2}}}\,\mathrm{S}(x)) + 
(1 + 2\,x)(4\,x^{4} + 9\,B\,x^{3}} \\
 & & \mbox{} + 9\,A\,x^{3} - 8\,C\,x^{3} + 3\,x^{3} - 12\,C\,x^{2
} + 9\,B\,x^{2} + 9\,A\,x^{2} + 3\,x^{2} - x \\
 & & \mbox{} - 6\,C\,x + 9\,A\,x + 9\,B\,x - C)({\frac {\partial 
}{\partial x}}\,\mathrm{S}(x))\mbox{} + 9\,(x - 1)^{2}\,B\,A\,
\mathrm{S}(x)=0\mbox{\hspace{14pt}}
\end{eqnarray*}
\end{maplelatex}

\end{maplegroup}
\begin{maplegroup}
\begin{mapleinput}
\mapleinline{active}{1d}{DE2:=sumdiffeq(second,k,S(x));}{%
}
\end{mapleinput}

\mapleresult
\begin{maplelatex}
\begin{eqnarray*}
\lefteqn{\mathit{DE2} := x\,(x - 1)\,(1 + x + x^{2})\,(1 + 2\,x)
^{2}\,({\frac {\partial ^{2}}{\partial x^{2}}}\,\mathrm{S}(x)) + 
(1 + 2\,x)(2\,x^{4} + 6\,b\,x^{4}} \\
 & & \mbox{} - 4\,d\,x^{4} + 6\,a\,x^{4} + 3\,b\,x^{3} + 3\,a\,x
^{3} + x^{3} - 6\,c\,x + 4\,x + 4\,d\,x + 2 \\
 & & \mbox{} - 3\,c)({\frac {\partial }{\partial x}}\,\mathrm{S}(
x))\mbox{} + ( - 12\,x^{4}\,a\,d + 36\,b\,a\,x^{4} + 4\,d^{2}\,x
^{4} - 12\,x^{4}\,b\,d \\
 & & \mbox{} - 2\,x^{3}\,d - 6\,x^{3}\,a\,d + 36\,b\,a\,x^{3} - 6
\,x^{3}\,b\,d + 9\,b\,a\,x^{2} - 12\,d\,x - 4\,d^{2}\,x \\
 & & \mbox{} + 12\,c\,x\,d + 6\,c\,d - 4\,d)\mathrm{S}(x)=0
\end{eqnarray*}
\end{maplelatex}

\end{maplegroup}
\begin{maplegroup}
\begin{mapleinput}

\mapleinline{active}{1d}{DE:=collect(collect(lhs(DE1)-lhs(DE2),S(x)),diff(S(x),x));
}{%
}
\end{mapleinput}

\mapleresult
\begin{maplelatex}
\begin{eqnarray*}
\lefteqn{\mathit{DE} := ((1 + 2\,x)(4\,x^{4} + 9\,B\,x^{3} + 9\,A
\,x^{3} - 8\,C\,x^{3} + 3\,x^{3} - 12\,C\,x^{2} + 9\,B\,x^{2}} \\
 & & \mbox{} + 9\,A\,x^{2} + 3\,x^{2} - x - 6\,C\,x + 9\,A\,x + 9
\,B\,x - C)\mbox{} - (1 + 2\,x)(2\,x^{4} \\
 & & \mbox{} + 6\,b\,x^{4} - 4\,d\,x^{4} + 6\,a\,x^{4} + 3\,b\,x
^{3} + 3\,a\,x^{3} + x^{3} - 6\,c\,x + 4\,x + 4\,d\,x \\
 & & \mbox{} + 2 - 3\,c))({\frac {\partial }{\partial x}}\,
\mathrm{S}(x))\mbox{} + (9\,(x - 1)^{2}\,B\,A + 12\,x^{4}\,a\,d
 - 36\,b\,a\,x^{4} \\
 & & \mbox{} - 4\,d^{2}\,x^{4} + 12\,x^{4}\,b\,d + 2\,x^{3}\,d + 
6\,x^{3}\,a\,d - 36\,b\,a\,x^{3} + 6\,x^{3}\,b\,d \\
 & & \mbox{} - 9\,b\,a\,x^{2} + 12\,d\,x + 4\,d^{2}\,x - 12\,c\,x
\,d - 6\,c\,d + 4\,d)\mathrm{S}(x)
\end{eqnarray*}
\end{maplelatex}

\end{maplegroup}
\begin{maplegroup}
\begin{mapleinput}

\mapleinline{active}{1d}{firstcoeff:=collect(frontend(coeff,[DE,S(x)]),x);
}{%
}
\end{mapleinput}

\mapleresult
\begin{maplelatex}
\begin{eqnarray*}
\lefteqn{\mathit{firstcoeff} := (12\,a\,d - 36\,b\,a - 4\,d^{2}
 + 12\,b\,d)\,x^{4}} \\
 & & \mbox{} + (2\,d + 6\,a\,d - 36\,b\,a + 6\,b\,d)\,x^{3} + (
 - 9\,b\,a + 9\,B\,A)\,x^{2} \\
 & & \mbox{} + ( - 18\,B\,A + 12\,d + 4\,d^{2} - 12\,c\,d)\,x + 9
\,B\,A - 6\,c\,d + 4\,d
\end{eqnarray*}
\end{maplelatex}

\end{maplegroup}
\begin{maplegroup}
\begin{mapleinput}

\mapleinline{active}{1d}{secondcoeff:=collect(frontend(coeff,[DE,diff(S(x),x)]),x);
}{%
}
\end{mapleinput}

\mapleresult
\begin{maplelatex}
\begin{eqnarray*}
\lefteqn{\mathit{secondcoeff} := (4 - 12\,b + 8\,d - 12\,a)\,x^{5
}} \\
 & & \mbox{} + (6 + 18\,B + 18\,A - 16\,C + 4\,d - 12\,a - 12\,b)
\,x^{4} \\
 & & \mbox{} + ( - 3\,b - 3\,a + 8 + 27\,B + 27\,A - 32\,C)\,x^{3
} \\
 & & \mbox{} + (12\,c - 7 - 8\,d - 24\,C + 27\,B + 27\,A)\,x^{2}
 \\
 & & \mbox{} + ( - 9 - 8\,C + 9\,A + 9\,B - 4\,d + 12\,c)\,x - C
 - 2 + 3\,c
\end{eqnarray*}
\end{maplelatex}

\end{maplegroup}
\begin{maplegroup}
\begin{mapleinput}
\mapleinline{active}{1d}{LIST:=\{coeffs(firstcoeff,x)\} union
\{coeffs(secondcoeff,x)\};}{%
}
\end{mapleinput}

\mapleresult
\begin{maplelatex}
\begin{eqnarray*}
\lefteqn{\mathit{LIST} := \{12\,a\,d - 36\,b\,a - 4\,d^{2} + 12\,
b\,d, \,9\,B\,A - 6\,c\,d + 4\,d, } \\
 & &  - 9\,b\,a + 9\,B\,A, \, - 18\,B\,A + 12\,d + 4\,d^{2} - 12
\,c\,d,  \\
 & & 6 + 18\,B + 18\,A - 16\,C + 4\,d - 12\,a - 12\,b, \,4 - 12\,
b + 8\,d - 12\,a,  \\
 & &  - 3\,b - 3\,a + 8 + 27\,B + 27\,A - 32\,C,  \\
 & &  - 9 - 8\,C + 9\,A + 9\,B - 4\,d + 12\,c,  \\
 & & 12\,c - 7 - 8\,d - 24\,C + 27\,B + 27\,A, \, - C - 2 + 3\,c
,  \\
 & & 2\,d + 6\,a\,d - 36\,b\,a + 6\,b\,d\}
\end{eqnarray*}
\end{maplelatex}

\end{maplegroup}
\begin{maplegroup}
\begin{mapleinput}
\mapleinline{active}{1d}{solve(LIST,\{A,B,C,a,b,c,d\});}{%
}
\end{mapleinput}

\mapleresult
\begin{maplelatex}
\begin{eqnarray*}
\lefteqn{\{d=d, \,B={\displaystyle \frac {1}{3}} \,d, \,A=
{\displaystyle \frac {1}{3}} \,d + {\displaystyle \frac {1}{3}} 
, \,b={\displaystyle \frac {1}{3}} \,d, \,C={\displaystyle 
\frac {1}{2}}  + {\displaystyle \frac {1}{2}} \,d, \,a=
{\displaystyle \frac {1}{3}} \,d + {\displaystyle \frac {1}{3}} 
, \,c={\displaystyle \frac {5}{6}}  + {\displaystyle \frac {1}{6}
} \,d\}, } \\
 & & \{d=d, \,B={\displaystyle \frac {1}{3}} \,d, \,A=
{\displaystyle \frac {1}{3}} \,d + {\displaystyle \frac {1}{3}} 
, \,b={\displaystyle \frac {1}{3}} \,d + {\displaystyle \frac {1
}{3}} , \,C={\displaystyle \frac {1}{2}}  + {\displaystyle 
\frac {1}{2}} \,d, \,a={\displaystyle \frac {1}{3}} \,d,  \\
 & & c={\displaystyle \frac {5}{6}}  + {\displaystyle \frac {1}{6
}} \,d\}, \{d=d, \,B={\displaystyle \frac {1}{3}} \,d + 
{\displaystyle \frac {1}{3}} , \,b={\displaystyle \frac {1}{3}} 
\,d, \,C={\displaystyle \frac {1}{2}}  + {\displaystyle \frac {1
}{2}} \,d, \,a={\displaystyle \frac {1}{3}} \,d + {\displaystyle 
\frac {1}{3}} ,  \\
 & & c={\displaystyle \frac {5}{6}}  + {\displaystyle \frac {1}{6
}} \,d, \,A={\displaystyle \frac {1}{3}} \,d\}, \{d=d, \,B=
{\displaystyle \frac {1}{3}} \,d + {\displaystyle \frac {1}{3}} 
, \,b={\displaystyle \frac {1}{3}} \,d + {\displaystyle \frac {1
}{3}} , \,C={\displaystyle \frac {1}{2}}  + {\displaystyle 
\frac {1}{2}} \,d,  \\
 & & a={\displaystyle \frac {1}{3}} \,d, \,c={\displaystyle 
\frac {5}{6}}  + {\displaystyle \frac {1}{6}} \,d, \,A=
{\displaystyle \frac {1}{3}} \,d\}
\;.
\mbox{\hspace{181pt}}
\end{eqnarray*}
\end{maplelatex}

\end{maplegroup}
\noindent
This leads to the unique solution
\[
\hypergeom{2}{1}{d/3,(1+d)/3}{(1+d)/2}
{1-\left(\frac{1-x}{1+2x}\right)^3}
=
(1+2x)^d\;\hypergeom{2}{1}{d/3,(1+d)/3}{(5+d)/6}{x^3}
\]
since the hypergeometric functions are symmetric w.r.t.\ their upper
parameters. Did you notice that this is exactly the computation from
\S~\ref{subsection:Polynomial Systems}? To solve this question, a nonlinear 
system had to be solved.

There is a theory of {\sl basic} hypergeometric ($q$-hypergeometric)
terms $A_k$ for which
$A_{k+1}/A_k$ is rational w.r.t.\ $q^k$. For almost all the results
and algorithms corresponding $q$-versions exist.

The corresponding series is called the {\sl basic hypergeometric series}
\[
_{r}\phi_{s}\left.\left(\begin{array}{cccc}
a_{1}&a_{2}&\cdots&a_{r}\\
b_{1}&b_{2}&\cdots&b_{s}\\
            \end{array}\right|q, x\right)
:=
\sum_{k=0}^\infty A_k\,x^{k}=
\]
\[
\sum_{k=0}^\infty \!\frac
{(a_{1};q)_{k}\cdot(a_{2};q)_{k}\cdots(a_{r};q)_{k}\,x^k}
{(b_{1};q)_{k}(b_{2};q)_{k}\cdots(b_{s};q)_{k}\qfac{q}{k}}
\left(\!(\!-1\!)^k q^{{k\choose 2}}\right)\!^{1+s-r}
\]
where $(a;q)_{k}:=\prod\limits_{j=0}^{k-1}{\left(1-a\,q^j\right)}$
is the {\sl $q$-Pochhammer-Symbol}.
$A_k$ is a
{\sl $q$-hypergeometric term} and fulfils the
recurrence equation $(k\in\N)$
\[
A_{k+1}:=\frac{(1-a_1q^k)\cdots (1-a_rq^k)}{(1-b_1q^k)\cdots
(1-b_sq^k)(1-q^{k+1})}\cdot A_k
\]
with the initial value $A_0:=1.$

All classical orthogonal polynomial families have at least one,
most families possess several $q$-analogues.

By the $q$-analogue of Zeilberger's algorithm, we get e.g.\ for the 
{\sl $q$-Laguerre polynomials}
\[
L_n^{(\alpha)}(x;q)=\frac{\qfac{q^{\alpha+1}}{n}}{\qfac{q}{n}}\,
\qphihyp{1}{1}{q^{-n}}{q^{\alpha+1}}{-xq^{n+\alpha+1}}
\]
the recurrence equation

\begin{maplegroup}

\begin{mapleinput}
\mapleinline{active}{1d}{read `qsum.mpl`;}{%
}
\end{mapleinput}

\mapleresult
\begin{maplelatex}
\[
\mathit{Copyright\ 1998,\ \ Harald\ Boeing\ \&\ Wolfram\ Koepf}
\]
\end{maplelatex}

\begin{maplelatex}
\[
\mathit{Konrad-Zuse-Zentrum\ Berlin}
\]
\end{maplelatex}

\end{maplegroup}
\begin{maplegroup}
\begin{mapleinput}

\mapleinline{active}{1d}{qsumrecursion(qpochhammer(q^(alpha+1),q,n)/qpochhammer(q,q,n)*}{%
}
\mapleinline{active}{1d}{qphihyperterm([q^(-n)],[q^(alpha+1)],q,-x*q^(n+alpha+1),k),}{%
}
\mapleinline{active}{1d}{q,k,L(n));
}{%
}
\end{mapleinput}

\mapleresult
\begin{maplelatex}
\begin{eqnarray*}
\lefteqn{ - q\,( - 1 + q^{n})\,\mathrm{L}(n) + ( - q^{2} + q^{(2
\,n + \alpha )}\,x - q + q^{(n + \alpha  + 1)} + q^{(1 + n)})\,
\mathrm{L}(n - 1)} \\
 & & \mbox{} + q\,(q - q^{(\alpha  + n)})\,\mathrm{L}(n - 2)=0
\;.
\mbox{\hspace{174pt}}
\end{eqnarray*}
\end{maplelatex}

\end{maplegroup}
\noindent
The $q$-analogue of Zeilberger's algorithm generates a third order
recurrence equation for the left-hand side of Jackson's $q$-analogue of
Dixon's identity
\[
\sum_{k=-n}^n
(-1)^k\,\qbinomial{n+b}{n+k}\,\qbinomial{n+c}{c+k}\,\qbinomial{b+c}{b+k}\,
q^{\frac{k(3k-1)}{2}}
=
\frac{
                        \qfac{q}{n+b+c}
                }{
                        \qfac{q}{n}\,\qfac{q}{b}\,\qfac{q}{c}
                }
\;:
\]

\begin{maplegroup}

\begin{mapleinput}

\mapleinline{active}{1d}{term:=(-1)^k*qbinomial(n+b,n+k,q)*qbinomial(n+c,c+k,q)*}{%
}
\mapleinline{active}{1d}{qbinomial(b+c,b+k,q)*q^(k*(3*k-1)/2);
}{%
}
\end{mapleinput}

\mapleresult
\begin{maplelatex}
\begin{eqnarray*}
\lefteqn{\mathit{term} := (-1)^{k}\,\mathrm{qbinomial}(n + b, \,n
 + k, \,q)\,\mathrm{qbinomial}(n + c, \,c + k, \,q)} \\
 & & \mathrm{qbinomial}(b + c, \,b + k, \,q)\,q^{(1/2\,k\,(3\,k
 - 1))}\mbox{\hspace{83pt}}
\end{eqnarray*}
\end{maplelatex}

\end{maplegroup}
\begin{maplegroup}
\begin{mapleinput}
\mapleinline{active}{1d}{RE:=qsumrecursion(term,q,k,S(n));}{%
}
\end{mapleinput}

\mapleresult
\begin{maplelatex}
\begin{eqnarray*}
\lefteqn{\mathit{RE} :=  - ( - q^{(2\,n)} + q)\,(q^{n} + 1)\,( - 
1 + q^{n})\,q^{3}\,\mathrm{S}(n) - (q^{5} - q^{(4 + n + c + b)}
 + q^{4}} \\
 & & \mbox{} + q^{(3 + 2\,n + b)} + q^{(3 + 2\,n + c)} - q^{(3 + 
2\,n)} - q^{(3 + c + b + n)} \\
 & & \mbox{} - q^{(3 + n + b)} - q^{(3 + n + c)} + q^{3} + q^{(2
\,n + b + 2)} + q^{(2\,n + c + 2)} \\
 & & \mbox{} - q^{(2\,n + c + 2 + b)} - q^{(3\,n + c + 2)} - q^{(
3\,n + 2)} - q^{(3\,n + b + 2)} \\
 & & \mbox{} + q^{(b + c + 2 + 4\,n)} - q^{(3\,n + 1)} + q^{(1 + 
4\,n + c + b)} + q^{(4\,n + c + b)})q \\
 & & \mathrm{S}(n - 1)\mbox{} + (q^{(2 + 2\,n)} - q^{(3\,n + c + 
b)} + q^{(2\,n + b + 2)} + q^{4} \\
 & & \mbox{} - q^{(4 + n + c + b)} - q^{(3\,n + 1 + c + b)} - q^{
(3 + n + b)} + q^{(3 + 2\,n + b)} + q^{6} \\
 & & \mbox{} + q^{5} - q^{(3\,n + c + 2 + b)} - q^{(4 + n + b)}
 + q^{(2\,n + c + 2)} - q^{(3 + n + c)} \\
 & & \mbox{} + q^{(3 + 2\,n + c)} - q^{(4 + n + c)})(q - q^{(n + 
c + b)})\,\mathrm{S}(n - 2)\mbox{} -  \\
 & & (q - q^{(n + c + b)})\,(q^{2} - q^{(n + b)})\,(q^{2} - q^{(n
 + c)})\,(q^{2} - q^{(n + c + b)}) \\
 & & \mathrm{S}(n - 3)=0
\;.
\end{eqnarray*}
\end{maplelatex}

\end{maplegroup}
\noindent
The $q$-analogue of Petkov\v{s}ek's algorithm decides whether a
$q$-holonomic recurrence
equation has $q$-hypergeometric term solutions.

It finds the right-hand side of the $q$-analogue of Dixon's identity:

\begin{maplegroup}

\begin{mapleinput}
\mapleinline{active}{1d}{qrecsolve(RE,q,S(n));}{%
}
\end{mapleinput}

\mapleresult
\begin{maplelatex}
\[
[[(q^{(n + 1)} - 1)\,\mathrm{S}(n + 1) + ( - q^{(1 + c + b + n)}
 + 1)\,\mathrm{S}(n)=0]]
\;.
\]
\end{maplelatex}

\end{maplegroup}
\noindent
In many cases,
much simpler is Paule's {\sl creative symmetrizing}. With this method,
we symmetrize the summand, and get the result in one step.

\begin{maplegroup}

\begin{mapleinput}

\mapleinline{active}{1d}{M:=qsimpcomb(subs(k=-k,term)/term,assume=[k,integer]);
}{%
}
\end{mapleinput}

\mapleresult
\begin{maplelatex}
\[
M := q^{k}
\]
\end{maplelatex}

\end{maplegroup}
\begin{maplegroup}
\begin{mapleinput}
\mapleinline{active}{1d}{qsumrecursion((1+M)/2*term,q,k,S(n));}{%
}
\end{mapleinput}

\mapleresult
\begin{maplelatex}
\[
(1 - q^{n})\,\mathrm{S}(n) + (q^{(n + c + b)} - 1)\,\mathrm{S}(n
 - 1)=0
\;.
\]
\end{maplelatex}

\end{maplegroup}
\noindent
There exist similar algorithms for definite integration instead of summation.

The {\sl Bateman integral representation}
\[
\int\limits_0^1 t^{c-1}\,(1-t)^{d-1}\,\hypergeom{2}{1}{a,b}{c}{tx}\,dt
=\frac{\Gamma(c)\Gamma(d)}{\Gamma(c+d)}\;
\hypergeom{2}{1}{a,b}{c+d}{x}
\]
is proven by

\begin{maplegroup}

\begin{mapleinput}

\mapleinline{active}{1d}{intrecursion(t^(c-1)*(1-t)^(d-1)*}{%
}
\mapleinline{active}{1d}{hyperterm([a,b],[c],t*x,k),t,B(k));
}{%
}
\end{mapleinput}

\mapleresult
\begin{maplelatex}
\[
 - (k + 1)\,(k + d + c)\,\mathrm{B}(k + 1) + \mathrm{B}(k)\,x\,(b
 + k)\,(a + k)=0
\]
\end{maplelatex}

\end{maplegroup}
\begin{maplegroup}
\begin{mapleinput}
\mapleinline{active}{1d}{assume(d>0,c>0);}{%
}
\mapleinline{active}{1d}{init:=int(t^(c-1)*(1-t)^(d-1),t=0..1);
}{%
}
\end{mapleinput}

\mapleresult
\begin{maplelatex}
\[
\mathit{init} := B(\mathit{c\symbol{126}}, \,\mathit{d
\symbol{126}})
\;.
\]
\end{maplelatex}

\end{maplegroup}

\noindent
We give our last example:
On top of the $q$-Askey-Wilson scheme \cite{KS}
we have the {\sl Askey-Wilson polynomials}
\[
p_n(x;a,b,c,d|q) = \qfac{ab}{n}\qfac{ac}{n}\qfac{ad}{n}\,a^{-n}
\cdot\; 
\]
\[
\qphihyp{4}{3}{q^{-n},abcdq^{n-1},ae^{i\,\theta},
ae^{-i\,\theta}}{a\,b,a\,c,a\,d}{q}
\;.
\]
The connection between those families can be written as
\[
p_n(x;\alpha,\beta,\gamma,d|q)
        \;=\;
\sum_{m=0}^n{C_m(n)\,p_m(x;a,b,c,d|q)}
\;,
\]
with {\sl connection coefficients} $C_m(n)$.

Askey and Wilson showed the following representation for $C_m(n)$
\[
C_m(n)=
\frac{\qfac{\alpha\,d,\beta\,d,\gamma\,d,q}{n}\,
\qfac{\alpha\,\beta\,\gamma\,d\,q^{n-1}}{m}}
{\qfac{\alpha\,d,\beta\,d,\gamma\,d,q,a\,b\,c\,d\,q^{m-1}}{m}\qfac{q}{n-m}}\,
q^{m^2-n\,m}\,d^{m-n}
\]
\[
\cdot\;
                \qphihyp{5}{4}{q^{m-n},\alpha\,\beta\,\gamma\,d\,q^{n+m-1},
                a\,d\,q^m,b\,d\,q^m,c\,d\,q^m}{a\,b\,c\,d\,q^{2\,m},
                \alpha\,d\,q^m,\beta\,d\,q^m,\gamma\,d\,q^m}{q}
\;.
\]
In the general case this is {\sl not} a $q$-hypergeometric term.

However,
for the special case $\beta=b$ and $\gamma=c$ we get by the $q$-Zeilberger
algorithm (\cite{BK}, compare \cite{GR}, Sections~7.5 and 7.6)
\[
                C_m(n)=
                \frac{
                        \qfac{\alpha/a}{n-m}
                        \qfac{\alpha b c d q^{n-1}}{m}
                        \qfac{q,b c,b d,c d}{n} a^{n-m}
                }{
                        \qfac{q,b c,b d,c d}{m}
                        \qfac{a b c d q^{m-1}}{m}
                        \qfac{q,a b c d q^{2 m}}{n-m}
                }
\;:
\]

\begin{maplegroup}

\begin{mapleinput}

\mapleinline{active}{1d}{c:='c': d:='d':}{%
}
\mapleinline{active}{1d}{term:=qpochhammer(alpha*d,beta*d,gamma*d,q,q,n)*}{%
}
\mapleinline{active}{1d}{qpochhammer(alpha*beta*gamma*d*q^(n-1),q,m)/}{%
}
\mapleinline{active}{1d}{qpochhammer(alpha*d,beta*d,gamma*d,q,a*b*c*d*q^(m-1),q,m)/}{%
}
\mapleinline{active}{1d}{qpochhammer(q,q,n-m)*q^(m^2-n*m)*d^(m-n)*}{%
}
\mapleinline{active}{1d}{qphihyperterm([q^(m-n),alpha*beta*gamma*d*q^(n+m-1),}{%
}
\mapleinline{active}{1d}{a*d*q^m,b*d*q^m,c*d*q^m],}{%
}
\mapleinline{active}{1d}{[a*b*c*d*q^(2*m),alpha*d*q^m,beta*d*q^m,gamma*d*q^m],q,q,j):}{%
}
\end{mapleinput}

\end{maplegroup}
\begin{maplegroup}
\begin{mapleinput}
\mapleinline{active}{1d}{qsumrecursion(subs(\{beta=b,gamma=c
\},term),q,j,C(m));}{%
}
\end{mapleinput}

\mapleresult
\begin{maplelatex}
\begin{eqnarray*}
\lefteqn{ - ( - 1 + q^{m})\,( - q + c\,d\,q^{m})\,(c\,a\,q^{(n + 
m)}\,b\,d - q)\,(b\,c\,q^{m} - q)\,( - q + b\,d\,q^{m})} \\
 & & ( - q^{3} + a\,b\,c\,d\,q^{(2\,m)})\,(a\,q^{m} - \alpha \,q
^{n})\,\mathrm{C}(m)\mbox{} + ( - q^{2} + c\,a\,q^{m}\,b\,d) \\
 & & ( - q + a\,b\,c\,d\,q^{(2\,m)})\,q^{2}\,( - q^{m} + q^{(1 + 
n)})\,( - q^{2} + \alpha \,c\,q^{(n + m)}\,b\,d) \\
 & & \mathrm{C}(m - 1)=0\mbox{\hspace{244pt}}
\end{eqnarray*}
\end{maplelatex}

\end{maplegroup}
\noindent
and similar results for $\alpha=a, \gamma=c$ and for
$\alpha=a, \beta=b$:

From these results one can derive the connection coefficients between many
families of the $q$-Askey-Wilson tableau by limit computations.


\begin{thebibliography}{10}

\bibitem{Beke1}
Beke, E.: Die Irreducibilit\"at der homogenen linearen
Differentialgleichungen. Math.\ Ann.\ 45, 1894, 278--294.

\bibitem{Beke2}
Beke, E.: Die symmetrischen Functionen bei linearen homogenen
Differentialgleichungen. Math.\ Ann.\ 45, 1894, 295--300.

\bibitem{BK}
B\"oing, H. and Koepf, W.: 
Algorithms for $q$-hypergeometric Summation in Computer Algebra.
J. Symbolic Computation, 1999, to appear.

\bibitem{BS}
Bronshtein, I.N. and Semedyayev, K.A.,
{\sl Handbook of Mathematics},
Springer, third edition 1985.

\bibitem{Maple}
Char, B.W.\ et al.: {\sl Maple V Language Reference Manual}.
Springer, New York, 1991.

\bibitem{FT}
Fleischer, J. and Tarasov, O.\ V.: Calculation of Feynman diagrams from
their small momentum expansion, Z.\ Phys.\ C64, 1994, 413.

\bibitem{GR}
Gasper, G.\ and Rahman, M.:
{\sl Basic Hypergeometric Series}.
Encyclopedia of Mathematics and its Applications,
35, Cambridge University Press, London and New York, 1990.

\bibitem{GCL}
Geddes, K. O., Czapor, S. R.\ and Labahn, G.:
{\sl Algorithms for Computer Algebra}. Kluwer Academic Publ.,
Boston/Dordrecht/London, 1992.

\bibitem{Gosper}
Gosper Jr., R. W.:
Decision procedure for indefinite hypergeometric
summation. Proc.\ Natl.\ Acad.\ Sci.\ USA 75, 1978, 40--42.

\bibitem{GK}
Gruntz, D. and Koepf, W.:
Maple package on formal power series.
Maple Technical Newsletter 2 (2), 1995, 22--28.

\bibitem{Reduce}
Hearn, A. C.:.
Reduce User's Manual, Version 3.6. RAND Co., Santa Monica, CA, 1995.

\bibitem{Axiom}
Jenks, R. D. and Sutor, R. S.: {\sl AXIOM. The Scientific Computation System}.
Springer, Berlin, 1993.

\bibitem{KS}
Koekoek, R.\ and Swarttouw, R. F.:
The Askey-scheme of hypergeometric orthogonal
polynomials and its $q$ analogue. Report \mbox{94--05,} Technische Universiteit
Delft, Faculty of Technical Mathematics and Informatics, Delft, 1994;
updated electronic version available at 
{\tt http://aw.twi.tudelft.nl/\verb+~+koekoek/research.html}.

\bibitem{Koepf1990}
Koepf, W.: 
Power series in Computer Algebra. J. Symbolic Computation
13, 1992, 581--603.

\bibitem{Koepf}
Koepf, W.: {\sl Hypergeometric Summation.
{\small\textsl 
An Algorithmic Approach to Summation and Special Function Identities.}}
Vieweg, Braunschweig/Wiesbaden, 1998.

\bibitem{Koornwinder}
Koornwinder, T. H.:
On Zeilberger's algorithm and its $q$-analogue: a rigorous description.
J.\ of Comput.\ and Appl.\ Math.\ 48, 1993, 91--111.

\bibitem{Macsyma}
{\sl Macsyma: Reference Manual}. Macsyma Inc., Arlington, MA 02174, USA.

\bibitem{MW}
Man, Y.-K.\ and Wright, F. J.:
Fast polynomial dispersion computation and its application to
indefinite summation.
Proc.\ of ISSAC 94, ACM Press, New York, 1994, 175--180.

\bibitem{SIAM}
Overhauser, A. W. and Kim, Y. I.:
Problem 94--2, SIAM Review 36, 1994, 107.

\bibitem{Derive}
Rich, A., Rich, J.\ and Stoutemyer, D.: {\sc Derive} {\sl User Manual},
Soft Warehouse, Inc.,
3660 Waialae Avenue, Suite 304, Honolulu, Hawaii, 96816-3236.

\bibitem{SZ}
Salvy, B.\ and Zimmermann, P.: GFUN: A package for the manipulation of
generating and holonomic functions in one variable. Rapports Techniques
143, INRIA, Rocquencourt (1992).

\bibitem{Sta}
Stanley, R.\ P.: Differentiably finite power series. Europ.\ J.\
Combinatorics 1, 1980, 175--188.

\bibitem{Mathematica}
Wolfram, St.: {\sl The Mathematica Book.} Wolfram Media, Champaign, Illinois,
and Cambridge University Press, Cambridge, 1996.


\end{thebibliography}
\end{document}